\documentclass[12pt]{amsart}
\usepackage{}

\usepackage{amsmath}
\usepackage{amsfonts}
\usepackage{amssymb}
\usepackage[all]{xy}           

\usepackage{bbding}
\usepackage{txfonts}
\usepackage{amscd}

\usepackage[shortlabels]{enumitem}
\usepackage{ifpdf}
\ifpdf
  \usepackage[colorlinks,final,
  hyperindex]{hyperref}
\else
  \usepackage[colorlinks,final,
  hyperindex]{hyperref}
\fi
\usepackage{tikz}
\usepackage[active]{srcltx}

\makeatletter

\newtheorem{thm}{Theorem}[section]
\newtheorem{lem}[thm]{Lemma}
\newtheorem{cor}[thm]{Corollary}
\newtheorem{pro}[thm]{Proposition}
\newtheorem{ex}[thm]{Example}
\newtheorem{rmk}[thm]{Remark}
\newtheorem{defi}[thm]{Definition}

\newcommand {\emptycomment}[1]{}

\newcommand{\kl}{\mathfrak l}
\newcommand{\kr}{\mathfrak r}

\newcommand{\bk}{\mathbb K}
\newcommand{\bc}{\mathbb C}

\newcommand{\Aut}{\mathrm{Aut}}

\newcommand{\End}{\mathrm{End}}

\newcommand{\Img}{\mathrm{Im}}
\newcommand{\Ker}{\mathrm{Ker}}

\newcommand{\HH}{\mathrm{HH}}

\def\id{\mathop {\fam0 id}\nolimits}
\def\tr{\triangleright}
\def\tl{\triangleleft}
\def\br{\blacktriangleright}
\def\bl{\blacktriangleleft}

\begin{document}

\title[Extending structures for perm algebras and perm bialgebras]
{Extending structures for perm algebras and perm bialgebras}

\author{Bo Hou}
\address{School of Mathematics and Statistics, Henan University, Kaifeng 475004,
China}
\email{bohou1981@163.com}
%
%


\begin{abstract}
We investigate the theory of extending structures by the unified product for
perm algebras, and the factorization problem as well as the classifying complements
problem in the setting of perm algebras. For a special extending structure,
non-abelian extension, we study the inducibility of a pair of automorphisms associated
to a non-abelian extension of perm algebras, and give the fundamental sequence of
Wells in the context of perm algebras. For a special extending structure,
bicrossed product, we introduce the concept of perm bialgebras, equivalently
characterized by Manin triples of perm algebras and certain matched pairs of perm algebras.
We introduce and study coboundary perm bialgebras, and our study leads to the
''$\mathcal{S}$-equation" in perm algebras, which is an analogue of the classical
Yang-Baxter equation. A symmetric solution of $\mathcal{S}$-equation gives a perm
bialgebra.
\end{abstract}


\maketitle

\tableofcontents 



\vspace{-4mm}

\section{Introduction}\label{sec:intr}

The notions of a perm algebra and perm operad have been introduced by Chapoton \cite{Cha}.
An associative algebra with the right-commutative identity is called perm
algebra. In \cite{MS} subalgebras of the perm algebra under the commutator
and the anti-commutator are considered. It turned out that every metabelian
Lie algebra can be embedded into perm algebra under the commutator. The perm operad
is the dual operad of the pre-Lie operad, which plays a specific role in the theory
of dialgebras \cite{Lod,GR}. Influenced by the results in the theory
of operads, Sartayev and his collaborators was considered a subalgebra of perm algebra
with derivation, and give a necessary and sufficient conditions for a left-symmetric
algebra can be embedded into perm algebra with derivation \cite{KS}; and proved that
every Novikov dialgebra can be embedded into perm algebra with derivation \cite{MS1}.
In \cite{Hur}, the author has considered the symmetric perm algebras.
In this paper, we study the extending structures problem of perm algebras and perm
bialgebras from the cohomological point of view.

Algebraic extension is an important method to study algebra structure. The extension problem,
induced for group by H\"{o}lder, is a famous and still open problem \cite{Hol}.
The extending structures problem was introduced and studied at the level of groups
in \cite{AM}. It was formulated for arbitrary categories in \cite{AM1}, where a partial
answer in the context of quantum groups was obtained. In recent years,
many researches have been conducted on extending structures problem of
various algebraic structures, such as Lie algebras \cite{AM3},
Leibinz algebras \cite{AM2}, associative algebras \cite{AM4}, left-symmetric
algebras \cite{Hon}, 3-Lie algebras \cite{Zha}, and so on. Here we study the
extending structures problem for perm algebras.

\vspace{2mm}
\noindent {\bf Extending structures problem.} {\it Let $A$ be a perm algebra
and $E$ a vector space containing $A$ as a subspace. Describe and classify up to an
isomorphism that stabilizes $A$ the set of all perm algebra structures
that can be defined on $E$ such that $A$ becomes a perm subalgebra of $E$.}

\vspace{2mm}
Equivalently, if we fix $V$ a complement of $A$ in the vector space $E$, the
extending structures problem can be restated by asking for the description and
classification of all perm algebras which contain and stabilize $A$ as a subalgebra
of codimension equal to ${\rm dim}(V)$. Therefore, we can study the extending
structures problem of perm algebras by introducing the unified product
and extending datum of perm algebras.

There is a close connection between the extending structures problem and many problems
related to the study of algebraic structure. First, if $A$ is the field $\bk$,
then the extending structures problem asks in fact for the classification of all
algebra structures on a given vector space $E$. In this case, the extending structures
problem is difficult. Commutative perm algebras are commutative associative algebras,
which is classified in \cite{BM,Bur} when the dimension is two. Here we prove
through direct calculation that two-dimensional noncommutative perm algebras have
only two, up to isomorphism.
Second, if the complement $V$ of $A$ in $E$ is one-dimensional, then the extension
algebra structure on $E$ is relatively easy. In this case, we can use some constants
and linear maps to classify perm algebra structure on $E$. Therefore, for a given perm
algebra $A$, we can specifically calculate all cases of extension algebra $E$.
Third, for some special extending datum, the extending algebra $E$ of $A$ is just
the non-abelian extension algebra of $A$, which is generalized abelian extension of $A$.
In this case, there is a close connection between the automorphism group of $A$ and
the automorphism group of $E$. Here, for a perm algebra $A$, based on the previous
classification results, we provide a necessary and sufficient condition that the
automorphism of $E$ can be obtained from the automorphisms of $A$ and $V$.
Fourth, the matched pairs of perm algebras are a special type of extending datum of
perm algebra, which is not only a powerful tool for solving the factorization problem,
but also closely related to the bialgebra structure. The factorization problem is an
old problem which stems in group theory \cite{Ore}. Here, we show that any perm
algebra which factorizes through two given perm algebras is isomorphic to a unified
product associated to a certain matched pair of perm algebras.

A bialgebra structure on a given algebraic structure is obtained as a coalgebra
structure together which gives the same algebraic structure on the dual space with
a set of compatibility. One of the most famous examples of bialgebras is the Lie
bialgebra \cite{Dri}. The antisymmetric infinitesimal bialgebras \cite{Agu, Bai1},
left-symmetric bialgebras \cite{Bai}, alternative D-bialgebras \cite{Gon}, Jordan
bialgebras \cite{Zhe}, 3-Lie bialgebras \cite{BCLM,BGS}, anti-flexible bialgebras
\cite{DBH}, have also been studied. In this paper, we give a bialgebra theory
for perm algebras. We show that the matched pairs of a perm algebra $A$ by the dual
space $A^{\ast}$ give a perm bialgebra structure on $A$, which is also an equivalent
structure of a Manin triple of perm algebras. Although we know very little about the
cohomology theory of perm algebras, we still consider the coboundary perm algebras,
following the research on left-symmetric bialgebras, because perm and
left-symmetric algebras have the same forms of dual bimodules. The study of such a class of
perm bialgebras also leads to the introduction of $\mathcal{S}$-equation in a perm algebra,
which is an analogue of the classical Yang-Baxter equation. A symmetric solution of
$\mathcal{S}$-equation in a perm algebra gives a perm bialgebra.

The paper is organized as follows. In Section \ref{sec: ext}, we recall the
notions of perm algebras and bimodules over perm algebras. We give the concept of an
extending datum of perm algebras with vector spaces and a unified product from an
extending datum, and give an answer to the extending structure problem for perm algebras.
We show that the bicrossed products of two perm algebras are special cases
of the unified product, which is used to solve the factorization problem for perm algebras.
Moreover, we consider the flag extending structure of perm algebras and give an answer to
the calculation of the classifying object of the flag extending structure of perm algebras.
In Section \ref{sec: non-abel}, we define the second non-abelian cohomology of perm algebras.
Applying the conclusion from the previous section, we can conclude that non-abelian extensions
of perm algebras are classified by the second non-abelian cohomology. In a given non-abelian
extension $\xymatrix@C=0.5cm{0\ar[r]& A\ar[r]^{\iota}& E\ar[r]^{\pi}& B\ar[r]&0}$ of perm
algebras, we show that the obstruction for a pair of automorphisms in
$\Aut(A)\times\Aut(B)$ to be induced by an automorphism in $\Aut(E)$ lies in the
second non-abelian cohomology. In particular, if the multiplication of $A$ is trivial,
we give the fundamental sequence of Wells in the context of perm algebras.
In Section \ref{sec: YB-equ}, we introduce the notion of perm bialgebras, and show that
perm bialgebras, standard Manin triples for perm algebras and certain matched pairs of
perm algebras are equivalent. The study of a special case of perm bialgebras,
the coboundary perm bialgebras, leads to the ''$\mathcal{S}$-equation" in perm algebras,
which is an analogue of the classical Yang-Baxter equation. A symmetric solution of
$\mathcal{S}$-equation gives a perm bialgebra.

Throughout this paper, we fix $\bk$ a field and characteristic zero. All the vector spaces
and algebras are over $\bk$, and all tensor products are also taking over $\bk$.

\section{Extending structures for perm algebras} \label{sec: ext}

In this section, we recall the notions of perm algebras and bimodules over perm algebras.
We introduce the definition of unified product for perm algebras, give a theoretical
answer to the extending structures problem and the classifying complements problem
for perm algebras.

\begin{defi}\label{def:anti}
Let $A$ be a $\bk$-vector space with a bilinear opreation $A\otimes A\rightarrow A$,
$a_{1}\otimes a_{2}\mapsto a_{1}a_{2}$ satisfying the (right commutative and associative)
identities:
$$
(a_{1}a_{2})a_{3}=a_{1}(a_{2}a_{3})=a_{1}(a_{3}a_{2}),
$$
for any $a_{1}, a_{2}, a_{3}\in A$. We call $A$ a {\rm perm algebra}.
\end{defi}

Here, the perm algebra is also called a right perm algebra. The left perm algebra is defined
by $a_{1}(a_{2}a_{3})=(a_{1}a_{2})a_{3}=(a_{2}a_{1})a_{3}$. Clearly, the opposite algebra
of a left (resp. right) perm algebra is a right (resp. left) perm algebra under the same
underlying vector space. Each commutative and associative algebra is a perm algebra.
And a perm algebra with identity is nothing but a unital commutative and associative algebra.
The perm algebra we consider in this paper does not require identity element.

\begin{ex}\label{ex:perm-alg}
$(i)$ Let $A$ be a commutative and associative algebra with a differential
$\mathbf{d}: A\rightarrow A$, i.e., $\mathbf{d}\circ\mathbf{d}=0$ and $\mathbf{d}(a_{1}a_{2})
=a_{1}\mathbf{d}(a_{2})+\mathbf{d}(a_{1})a_{2}$, for $a_{1}, a_{2}\in A$. Define
$\cdot: A\otimes A\rightarrow A$ by $a_{1}\otimes a_{2}\mapsto a_{1}\mathbf{d}(a_{2})$.
Then $(A, \cdot)$ is a perm algebra.

$(ii)$ Assume $V$ is a $\bk$ vector space and $S(V)$ is the free symmetric algebra on $V$.
Consider the vector space $R:=V\otimes S(V)$, and define bilinear opreation $R\otimes
R\rightarrow R$ by $(v_{1}, x_{1})(v_{2}, x_{2})=v_{1}\otimes(x_{1}v_{2}x_{2})$, for any
$(v_{1}, x_{1}), (v_{2}, x_{2})\in R$. Then $R$ is a perm algebra, is called the
free perm algebra on $V$.
\end{ex}

Let $A$ and $B$ be two perm algebras. A linear map $f: A\rightarrow B$ is called a
{\it morphism of perm algebras} if for any $a_{1}, a_{2}\in A$, $f(a_{1},a_{2})=
f(a_{1})f(a_{2})$. A morphism $f$ is said to be an isomorphism if $f$ is a bijection.
We denote the automorphism group of $A$ by $\Aut(A)$.

\begin{ex}\label{ex:perm-algs}
Let $A=\bc\{e_{1}, e_{2}\}$ be a 2-dimensional vector space over complex field $\bc$.
We consider the noncommutative perm algebra structure on $A$. Suppose
\begin{align*}
& e_{1}e_{1}=s_{1}e_{1}+s_{2}e_{2},\qquad\qquad e_{2}e_{2}=t_{1}e_{1}+t_{2}e_{2},\\
& e_{1}e_{2}=p_{1}e_{1}+p_{2}e_{2},\qquad\qquad e_{2}e_{1}=q_{1}e_{1}+q_{2}e_{2}.
\end{align*}
where $s_{1}, s_{2}, p_{1}, p_{2}, q_{1}, q_{2}, t_{1}, t_{2}\in\bc$.
Since $(xy)z=x(yz)=x(zy)$, for any $x, y, z\in A$, we get
\begin{align*}
&\qquad\quad p_{1}q_{2}=q_{1}q_{2}=t_{1}s_{2}=p_{1}p_{2},\qquad
p_{1}p_{1}+t_{1}p_{2}=s_{1}t_{1}+p_{1}t_{2},\\
&\qquad p_{1}q_{1}=q_{1}q_{1},\qquad t_{1}q_{2}=t_{1}p_{2},\qquad t_{2}p_{2}=t_{2}q_{2},
\qquad s_{1}p_{1}=s_{1}q_{1},\\
&s_{1}p_{2}+t_{2}s_{2}=p_{1}s_{2}+p_{2}p_{2}=q_{1}s_{2}+p_{2}q_{2},\qquad
q_{1}s_{2}+q_{2}q_{2}=s_{1}q_{2}+t_{2}s_{2}.
\end{align*}
Moreover, since the multiplication is noncommutative, we will discuss the following
situations.

$(1)$ $p_{1}\neq q_{1}$. From the above equations, it can be directly obtained that
$s_{1}=q_{1}=q_{2}=0$, $t_{2}s_{2}=t_{1}s_{2}=p_{1}p_{2}=t_{1}p_{2}=t_{2}p_{2}=0$,
$p_{1}s_{2}+p_{2}p_{2}=0$ and $p_{1}p_{1}+t_{1}p_{2}=p_{1}t_{2}$.
$(i)$ If $t_{2}=s_{2}=0$. We get $p_{1}=p_{2}=0$. This contradicts $p_{1}\neq q_{1}$.
$(ii)$ If $t_{2}=0$ and $s_{2}\neq0$. We get $t_{1}=p_{1}=0$.
This contradicts $p_{1}\neq q_{1}$. $(iii)$ If $s_{2}=0$ and $t_{2}\neq0$.
we get $p_{2}=0$ and $p_{1}p_{1}=p_{1}t_{2}$. Since $p_{1}\neq q_{1}=0$ in this case,
we have $p_{1}=t_{2}$. Thus, in this case, the perm multiplication on $A$ is given by
$$
e_{1}e_{1}=0,\qquad e_{1}e_{2}=t_{2}e_{1},\qquad e_{2}e_{1}=0,\qquad
e_{2}e_{2}=t_{1}e_{1}+t_{2}e_{2},
$$
where $t_{1}, t_{2}\in\bc$ and $t_{2}\neq0$.

$(2)$ $p_{2}\neq q_{2}$. The equations above are transformed into $t_{1}=t_{2}
=p_{1}=q_{1}=0$, $q_{2}q_{2}=s_{1}q_{2}$ and $s_{1}p_{2}=p_{2}p_{2}=p_{2}q_{2}$.
$(i)$ If $q_{2}=0$, we get $p_{2}=0$. This contradicts $p_{2}\neq q_{2}$.
$(ii)$ If $q_{2}=s_{1}\neq0$. We obtain $p_{2}p_{2}=p_{2}q_{2}$. Since $p_{2}\neq q_{2}$
in this case, we get $p_{2}=0$. Thus, in this case, the perm multiplication on $A$ is given by
$$
e_{1}e_{1}=s_{1}e_{1}+s_{2}e_{2},\qquad e_{1}e_{2}=0,\qquad e_{2}e_{1}=s_{1}e_{2},\qquad
e_{2}e_{2}=0,
$$
where $s_{1}, s_{2}\in\bc$ and $s_{1}\neq0$. Note that the perm algebras obtained
in $(1)$ are isomorphic to the perm algebras obtained in $(2)$, we get that,
up to isomorphism, the noncommutative perm algebras have only one class.
We denote this 2-dimensional perm algebra by $A_{s_{1}}$ if $s_{2}=0$, and by
$A_{s_{1}, s_{2}}$ if $s_{2}\neq0$. Clearly, $A_{s_{1}}$ is not isomorphic to
$A_{s'_{1}, s'_{2}}$, for any nonzero $s_{1}, s_{2}\in\bc$. For any nonzero
$s_{1}, s'_{1}\in\bc$, one can check that $f: A_{s_{1}}\rightarrow A_{s'_{1}}$,
$f(e_{1})=\frac{s_{1}}{s'_{1}}(e_{1}+e_{2})$, $f(e_{2})=e_{2}$, is an isomorphism.
For any nonzero $s_{1}, s_{2}, s'_{1}, s'_{2}\in\bc$, one can check that
$f: A_{s_{1}, s_{2}}\rightarrow A_{s'_{1}, s'_{2}}$,
$f(e_{1})=\frac{s_{1}}{s'_{1}}(e_{1}+e_{2})$, $f(e_{2})=\frac{s_{1}s_{1}s'_{2}}
{s'_{1}s'_{1}s_{2}}e_{2}$, is an isomorphism. This means, each noncommutative
2-dimensional perm $\bc$-algebra is isomorphic to one of the following perm algebra
on $\bc\{e_{1}, e_{2}\}$:
\begin{itemize}
\item[$(i)$] $e_{1}e_{1}=e_{1}$,\qquad $e_{1}e_{2}=0$,\qquad $e_{2}e_{1}=e_{2}$,\qquad
     $e_{2}e_{2}=0$;
\item[$(ii)$] $e_{1}e_{1}=e_{1}+e_{2}$,\qquad $e_{1}e_{2}=0$,\qquad $e_{2}e_{1}=e_{2}$,\qquad
     $e_{2}e_{2}=0$.
\end{itemize}
\end{ex}

Now we introduce the definition of bimodule over a perm algebra.

\begin{defi}\label{def:anti-mod}
Let $A$ be a perm algebra, $V$ be a $\bk$-vector space. If there exist a pair of
bilinear maps $l: A\otimes V\rightarrow V$, $l(a, v)=av$, and $r: V\otimes A
\rightarrow V$, $r(v, a)=va$, such that for any $a_{1}, a_{2}\in A$ and $v\in V$,
\begin{align*}
&\qquad v(a_{1}a_{2})=(v a_{1})a_{2}=(v a_{2})a_{1},\\
&(a_{1}a_{2})v=a_{1}(a_{2}v)=a_{1}(va_{2})=(a_{1}v)a_{2}.
\end{align*}
We call $(V, l, r)$ (or simply $V$) a {\rm bimodule} over $A$, or $V$ an $A$-bimodule.
\end{defi}

By direct calculations, we can give an equivalent condition as follows.

\begin{pro}\label{pro:rep}
Let $A$ be a perm algebra, $V$ be a $\bk$-vector space, $l: A\otimes V\rightarrow V$,
$l(a, v)=av$, and $r: V\otimes A\rightarrow V$, $r(v, a)=va$, be two bilinear maps.
Then $(V, l, r)$ is an $A$-bimodule if and only if $A\oplus V$ is a perm algebra
under the following operations:
$$
(a_{1}, v_{1})(a_{2}, v_{2})=\Big(a_{1}a_{2},\, a_{1}v_{2}+v_{1}a_{2}\Big),
$$
for all $a_{1}, a_{2}\in A$ and $v_{1}, v_{2}\in V$. This perm algebra is called the
{\rm semidirect product} of $A$ by bimodule $V$, denoted by $A\ltimes V$.
\end{pro}

Moreover, for convenience, we sometimes need to consider the representations of perm algebras.

\begin{defi}\label{def:anti-rep}
Let $A$ be a perm algebra, $V$ be a $\bk$-vector space, and $\kl, \kr: A\rightarrow
\End_{\bk}(V)$ be two linear maps. The triple $(V, \kl, \kr)$ is called a {\rm representation}
of $A$, if for any $a_{1}, a_{2}\in A$,
\begin{align*}
&\qquad\kr(a_{1}a_{2})=\kr(a_{2})\circ\kr(a_{1})=\kr(a_{1})\circ\kr(a_{2}),\\
&\kl(a_{1}a_{2})=\kl(a_{1})\circ\kl(a_{2})=\kl(a_{1})\circ\kr(a_{2})=\kr(a_{2})\circ\kl(a_{1}).
\end{align*}
\end{defi}

Let $(V, \kl, \kr)$ and $(V', \kl', \kr')$ be two representations of a perm algebra $A$.
A linear map $f: V\rightarrow V$ is called a {\it morphism of representations} if
for any $v\in V$, $f\circ\kl(v)=\kl'(v)\circ f$ and $f\circ\kr(v)=\kr'(v)\circ f$.
If the morphism $f$ is a bijection, we call that $f$ is an isomorphism and
the representations $(V, \kl, \kr)$ and $(V', \kl', \kr')$ are isomorphic.
Indeed, the notions of bimodule and representation of a perm algebra is equivalent.

\begin{pro}\label{pro:mod-rep}
Let $A$ be a perm algebra, $V$ be a $\bk$-vector space, and $\kl, \kr: A\rightarrow
\End_{\bk}(V)$ be two linear maps. Then $(V, \kl, \kr)$ is a representation of $A$ if
and only if $(V, l, r)$ is a bimodule over $A$, where $l: A\otimes V\rightarrow V$,
$l(a, v)=\kl(a)(v)$, and $r: V\otimes A\rightarrow V$, $r(v, a)=\kr(a)(v)$.

Moreover, a linear map $f$ is a morphism from representation $(V, \kl, \kr)$ to
representation $(V', \kl', \kr')$ if and only if $f$ is a morphism from bimodule $(V, l, r)$
to bimodule $(V', l', r')$.
\end{pro}

For a perm algebra $A$, we get linear maps $\kl_{A}: A\rightarrow\End_{\bk}(A)$,
$\kl_{A}(a_{1})(a_{2})=a_{1}a_{2}$, $\kr_{A}: A\rightarrow\End_{\bk}(A)$,
$\kr_{A}(a_{1})(a_{2})=a_{2}a_{1}$ for any $a_{1}, a_{2}\in A$ by the multiplication of $A$.
It is easy to see that $(A, \kl_{A}, \kr_{A})$ is a representation of $A$.
We call this representation the {\it regular representation} of $A$.

\subsection{Unified product for perm algebras} \label{subsec:unif}

\begin{defi}\label{Def:exten}
Let $A$ be a given perm algebra, $E$ be a vector space. An {\rm extending structure}
of $A$ by $V$ is a perm algebra $E$ such that $A$ is a subalgebra of $E$ and $V$ a
complement of $A$ in $E$, which fit into the following exact sequence as vector spaces:
$$
\xymatrix{ 0 \ar[r] & A \ar[r]^{\iota} & E \ar[r]^{\pi} & V \ar[r] & 0 }.
$$
\end{defi}

The extending problem is to describe and classify up to an isomorphism the set of all perm
algebra structures that can be defined on $E$ such that $A$ is a perm subalgebra of $E$.

\begin{defi} \label{Def:exte-equ}
Let $A$ be a perm algebra, $E$ be a vector space. Suppose $A$ is a subspace of $E$
and $V$ is a complement of $A$ in $E$. That is, $E\cong A\oplus V$ as vector spaces.
Consider the diagram in the category of vector spaces:
$$
\xymatrixrowsep{0.5cm}
\xymatrix{0\ar[r]& A\ar[r]^{\iota}
\ar@{=}[d]& E\ar[r]^{\pi}\ar[d]^{\varphi} & V\ar[r]\ar@{=}[d]&0\;\\
0\ar[r]& A\ar[r]^{\iota} & E\ar[r]^{\pi} & V\ar[r]&0,}
$$
where $\pi: E\rightarrow V$ is the canonical projection and $\iota: A\rightarrow E$ is the
inclusion map. We say that the linear map $\varphi: E\rightarrow E$ {\rm stabilizes}
$A$ if the left square of the diagram is commutative, and $\varphi: E\rightarrow E$
{\rm costabilizes} $V$ if the right square of the diagram is commutative.

Let $(E, \diamond)$ and $(E, \diamond')$ be two perm algebra structures
on $E$ both containing $A$ as a subalgebra. $(E, \diamond)$ and $(E, \diamond')$ are called
{\rm equivalent}, and denoted by $(E, \diamond)\approx(E, \diamond')$, if there exists a
perm algebra isomorphism $\varphi: (E, \diamond)\rightarrow(E, \diamond')$ which stabilizes
${Z}$. Denote by $\mathrm{Ext}(E, A)$ the set of equivalent classes of $A$ through $V$. Perm
algebras $(E, \diamond)$ and $(E, \diamond')$ are called {\rm cohomologous}, and we denote this
by $(E, \diamond)\equiv(E, \diamond')$, if there exists a perm algebra isomorphism $\varphi:
(E, \diamond)\rightarrow(E, \diamond')$ which stabilizes $A$ and costabilizes $V$. We denote
by $\mathbf{Ext}(E, A)$ the set of cohomologous classes of $A$ through $V$.
\end{defi}

In the following of this section we give a theoretical answer to the
extending problem by introducing the unified product for perm algebras.

\begin{defi}\label{Def:ext-datum}
Let $A$ be a perm algebra and $V$ be a $\bk$-vector space. An {\rm extending datum}
of $A$ by $V$ is a system $\Omega(A, V)=(\br, \bl, \tr, \tl, \chi, \cdot)$ consisting
bilinear linear maps $\br: A\times V\rightarrow V$, $\bl: V\times A\rightarrow V$,
$\tr: V\times A\rightarrow A$, $\tl: A\times V\rightarrow A$,
$\chi: V\times V\rightarrow A$, and $\cdot: V\times V\rightarrow V$.

Let $\Omega(A, V)=(\br, \bl, \tr, \tl, \chi, \cdot)$ be an extending datum.
Denote by $A\natural_{\Omega(A, V)}V$, or simply $A\natural V$, the vector space
$A\oplus V$ with the bilinear map $\ast: (A\oplus V)\times (A\oplus V)\rightarrow
A\oplus V$ defined by
\begin{equation}\label{Equ:prod}
(a_{1}, v_{1})\ast(a_{2}, v_{2})=\Big(a_{1}a_{2}+v_{1}\tr a_{2}+a_{1}\tl v_{2}
+\chi(v_{1}, v_{2}),\ \ v_{1}\cdot v_{2}+a_{1}\br v_{2}+v_{1}\bl a_{2}\Big),
\end{equation}
for all $a_{1}, a_{2}\in A$, $v_{1}, v_{2}\in V$. $A\natural V$ is called the {\rm unified
product} of $A$ and $V$ if $(A\oplus V,\; \ast)$ a perm algebra. In this case, the extending
datum $\Omega(A, V)$ is called an {\rm perm extending structure} of $A$ by $V$. We denote by
$\mathfrak{E}(A, V)$ the set of all perm extending structures of $A$ by $V$.
\end{defi}

According to this definition, we can directly give the equivalent characterization of
the unified product.

\begin{pro}\label{Pro:exten}
Let $A$ be a perm algebra, $V$ be a vector space and $\Omega(A, V)$ be an extending datum of
$A$ by $V$. Then, $A\natural V$ is a unified product if and only if the following conditions
hold for any $a, a_{1}, a_{2}\in A$, $v, v_{1}, v_{2}, v_{3}\in V$,
\begin{align}
(v\bl a_{1})\bl a_{2}&=v\bl(a_{1}a_{2})=v\bl(a_{2}a_{1}),                      \label{ext1}\\
(a_{1}a_{2})\br v&=a_{1}\br(a_{2}\br v)
=a_{1}\br(v\bl a_{2})=(a_{1}\br v)\bl a_{2},                                   \label{ext2}\\
v\tr(a_{1}a_{2})&=v\tr(a_{2}a_{1})=(v\tr a_{1})a_{2}+(v\bl a_{1})\tr a_{2},    \label{ext3}\\
(a_{1}a_{2})\tl v&=a_{1}(a_{2}\tl v)+a_{1}\tl(a_{2}\br v)                      \label{ext4}\\
&=a_{1}(v\tr a_{2})+a_{1}\tl(v\bl a_{2})
=(a_{1}\tl v)a_{2}+(a_{1}\br v)\tr a_{2},                                      \nonumber\\
a\br(v_{1}\cdot v_{2})&=a\br(v_{2}\cdot v_{1})
=(a\br v_{1})\cdot v_{2}+(a\tl v_{1})\br v_{2},                                \label{ext5}\\
(v_{1}\cdot v_{2})\bl a&=v_{1}\cdot(v_{2}\bl a)+v_{1}\bl(v_{2}\tr a)           \label{ext6}\\
&=v_{1}\cdot(a\br v_{2})+v_{1}\bl(a\tl v_{2})
=(v_{1}\bl a)\cdot v_{2}+(v_{1}\tr a)\br v_{2},                                \nonumber\\
a\chi(v_{1}, v_{2})+a\tl(v_{1}\cdot v_{2})&=a\chi(v_{2}, v_{1})
+a\tl(v_{2}\cdot v_{1})=(a\tl v_{1})\tl v_{2}+\chi(a\br v_{1}, v_{2}),       \label{ext7}
\end{align}
\begin{align}
\chi(v_{1}, v_{2})a+(v_{1}\cdot v_{2})\tr a
&=v_{1}\tr(v_{2}\tr a)+\chi(v_{1}, v_{2}\bl a)                                 \label{ext8}\\
&=v_{1}\tr(a\tl v_{2})+\chi(v_{1}, a\br v_{2})
=(v_{1}\tr a)\tl v_{2}+\chi(v_{1}\bl a, v_{2}),                                 \nonumber\\
(v_{1}\cdot v_{2})\cdot v_{3}+\chi(v_{1}, v_{2})\br v_{3}
&=v_{1}\cdot(v_{2}\cdot v_{3})+v_{1}\bl\chi(v_{2}, v_{3})
=v_{1}\cdot(v_{3}\cdot v_{2})+v_{1}\bl\chi(v_{3}, v_{2}),                      \label{ext9}\\
\chi(v_{1}, v_{2})\tl v_{3}+\chi(v_{1}\cdot v_{2}, v_{3})
&=v_{1}\tr\chi(v_{2}, v_{3})+\chi(v_{1}, v_{2}\cdot v_{3})
=v_{1}\tr\chi(v_{3}, v_{2})+\chi(v_{1}, v_{3}\cdot v_{2}).                     \label{ext10}
\end{align}
\end{pro}

\begin{proof}
Define
\begin{align*}
R_{1}((a_{1}, v_{1}), (a_{2}, v_{2}), (a_{3}, v_{3}))
&=(a_{1}, v_{1})\ast((a_{2}, v_{2})\ast(a_{3}, v_{3}))
-((a_{1}, v_{1})\ast(a_{2}, v_{2}))\ast(a_{3}, v_{3})\\
R_{2}((a_{1}, v_{1}), (a_{2}, v_{2}), (a_{3}, v_{3}))
&=(a_{1}, v_{1})\ast((a_{2}, v_{2})\ast(a_{3}, v_{3}))
-(a_{1}, v_{1})\ast((a_{3}, v_{3})\ast(a_{2}, v_{2})),
\end{align*}
for $(a_{1}, v_{1}), (a_{2}, v_{2}), (a_{3}, v_{3})\in A\oplus V$.
Then $A\natural V$ is a perm algebra if and only if
$R_{1}((a_{1}, v_{1})$, $(a_{2}, v_{2}), (a_{3}, v_{3}))
=R_{2}((a_{1}, v_{1})$, $(a_{2}, v_{2}), (a_{3}, v_{3}))=0$
for all $(a_{1}, v_{1}), (a_{2}, v_{2}), (a_{3}, v_{3})\in A\oplus V$.
By a direct computation, it is easy to see that,
$R_{1}((0, v), (a_{1}, 0), (a_{2}, 0))=R_{2}((0, v), (a_{1}, 0), (a_{2}, 0))=0$ if and only if
Eqs. (\ref{ext1}) and (\ref{ext3}) hold;
$R_{1}((a_{1}, 0), (0, v), (a_{2}, 0))=R_{2}((a_{1}, 0), (0, v), (a_{2}, 0))=0$ and
$R_{1}((a_{1}, 0), (a_{2}, 0), (0, v))=R_{2}((a_{1}, 0), (a_{2}, 0), (0, v))=0$, if and only if
Eqs. (\ref{ext2}) and (\ref{ext4}) hold;
$R_{1}((a, 0), (0, v_{1}), (0, v_{2}))=R_{2}((a, 0), (0, v_{1}), (0, v_{2}))=0$ if and only if
Eqs. (\ref{ext5}) and (\ref{ext7}) hold;
$R_{1}((0, v_{1}), (a, 0), (0, v_{2}))=R_{2}((0, v_{1}), (a, 0), (0, v_{2}))=0$ and
$R_{1}((0, v_{1}), (0, v_{2}), (a, 0))=((0, v_{1}), (0, v_{2}), (a, 0))=0$, if and only if
Eqs. (\ref{ext6}) and (\ref{ext8}) hold;
$R_{1}((0, v_{1}), (0, v_{2}), (0, v_{3}))=R_{2}((0, v_{1}), (0, v_{2}), (0, v_{3}))=0$
if and only if Eqs. (\ref{ext9}) and (\ref{ext10}) hold.
The proof is finished.
\end{proof}

It is easy to see that, if $A\natural V$ is a unified product, then $A$ is a subalgebras
of perm algebra $A\natural V$ and $(V, \br, \bl)$ is a bimodule over perm algebra $A$.
Conversely, we can also prove that any perm algebra structure on the vector space $E$
containing $A$ as a subalgebra is isomorphic to a unified product.

\begin{thm}\label{Thm:ext-unif}
Let $A$ be a perm algebra, $E$ be a vector space containing $A$ as a subspace and 
$``\diamond"$ a perm algebra structure on $E$ such that $A$ is a subalgebra of 
$(E, \diamond)$. Then, there exists a perm extending structure $\Omega(A, V)=(\br, \bl, 
\tr, \tl, \chi, \cdot)$ of $A$ by a subspace $V$ of $E$ and an isomorphism of perm 
algebras $(E, \diamond)\cong A\natural V$ which stabilizes $A$ and co-stabilizes $V$.
\end{thm}

\begin{proof}
Let $\pi_{A}: E\rightarrow A$ be the projection map such that $\pi_{A}(a)=a$ for all $a\in A$.
Set $V=\Ker(\pi_{A})$ which is a complement of $A$ in $E$. Then, we present an extending
datum $\Omega(A, V)=(\br, \bl$, $\tr, \tl, \chi, \cdot)$ of $A$ by a subspace
$V$ of $E$ defined as follows:
\begin{align*}
\br: A\times V \rightarrow V, &\qquad\quad a\br v=a\diamond v-\pi_{A}(a\diamond v),\\
\bl: V\times A \rightarrow V, &\qquad\quad v\bl a=v\diamond a-\pi_{A}(v\diamond a),\\
\tr: V\times A \rightarrow A, &\qquad\quad v\tr a=\pi_{A}(v\diamond a),\\
\tl: A\times V \rightarrow A, &\qquad\quad a\tl v=\pi_{A}(a\diamond v),\\
\chi: V\times V \rightarrow A, &\qquad\quad \chi(v_{1}, v_{2})=\pi_{A}(v_{1}\diamond v_{2}),\\
\cdot: V\times V \rightarrow V, &\qquad\quad v_{1}\cdot v_{2}=v_{1}\diamond v_{2}
-\pi_{A}(v_{1}\diamond v_{2}),
\end{align*}
for any $a\in A$, $v, v_{1}, v_{2}\in V$. It is easy to see that
$\varphi: A\oplus V\rightarrow E$, $\varphi(a, v)=a+v$ is a linear isomorphism,
whose inverse is $\varphi^{-1}(e)=(\pi_{A}(e), e-\pi_{A}(e))$ for all $e\in E$.
Note that the linear isomorphism $\varphi$ induces a unique perm product
on $A\oplus V$, which is given by
$$
(a_{1}, v_{1})(a_{2}, v_{2})=\varphi^{-1}\Big(\varphi(a_{1}, v_{1})\cdot
\varphi(a_{2}, v_{2})\Big)
$$
such that $\varphi$ is an isomorphism of perm algebras. To prove
$\Omega(A, V)=(\br, \bl, \tr, \tl, \chi, \cdot)$ is a perm extending
structure of $A$ by $V$, by Proposition \ref{Pro:exten}, we need to show that this
perm product is just the one defined by Eq. (\ref{Equ:prod}) associated
to $\Omega(A, V)$. Indeed, for all $(a_{1}, v_{1}), (a_{2}, v_{2})\in A\oplus V$, we have
\begin{align*}
&\; (a_{1}, v_{1})(a_{2}, v_{2})\\
=&\; \varphi^{-1}\Big(\varphi(a_{1}, v_{1})\diamond\varphi(a_{2}, v_{2})\Big)\\
=&\; \varphi^{-1}\Big(a_{1}a_{2}+a_{1}\diamond v_{2}+v_{1}\diamond a_{2}
+v_{1}\diamond v_{2}\Big)\\
=&\; \Big(a_{1}a_{2}+\pi_{A}(a_{1}\diamond v_{2})+\pi_{A}(v_{1}\diamond a_{2})
+\pi_{A}(v_{1}\diamond v_{2}),\\[-2mm]
&\qquad a_{1}\diamond v_{2}+v_{1}\diamond a_{2}+v_{1}\diamond v_{2}
-\pi_{A}(a_{1}\diamond v_{2})-\pi_{A}(v_{1}\diamond a_{2})-\pi_{A}(v_{1}\diamond v_{2})\Big)\\
=&\; \Big(a_{1}a_{2}+v_{1}\tr a_{2}+a_{1}\tl v_{2}+\chi(v_{1}, v_{2}),\ \
v_{1}\cdot v_{2}+a_{1}\br v_{2}+v_{1}\bl a_{2}\Big).
\end{align*}
Thus, $\varphi$ is an isomorphism of perm algebras. Moreover, we consider the diagram
$$
\xymatrixrowsep{0.5cm}
\xymatrix{0\ar[r]& A\ar[r]^{\iota_{A}\;}\ar@{=}[d]& A\natural
V\ar[r]^{\;\pi_{V}}\ar[d]^{\varphi}& V\ar[r]\ar@{=}[d]&0\\
0\ar[r]& A\ar[r]^{\iota} & E\ar[r]^{\pi} & V\ar[r]&0}
$$
in the category of vector spaces, where $\iota_{A}(a)=(a, 0)$, $\pi_{V}(a, v)=v$.
It is easy to see that the diagram is commutative. Thus, the proof is finished.
\end{proof}

For classifying all perm algebra structures on $E$ containing $A$ as a subalgebra, by
Theorem \ref{Thm:ext-unif}, we noly need to classify all unified products $A\natural V$
associated to perm extending structures $\Omega(A, V)=(\br, \bl, \tr, \tl, \chi, \cdot)$
for a given complement $V$ of $A$ in $E$. Next, we construct the cohomological objects
for extending structures which will parameterize the classifying sets $\mathrm{Ext}(E, A)$
and $\mathbf{Ext}(E, A)$, respectively.

\begin{defi} \label{Def:unif-coho}
Let $A$ be a perm algebra and $V$ be a vector space. If there exists a linear map $\lambda:
V\rightarrow A$ and a linear bijiction $\rho: V\rightarrow V$, such that the perm extending
structure $\Omega(A, V)=(\br, \bl, \tr, \tl, \chi, \cdot)$ can be obtained from another
extending structure $\Omega'(A, V)=(\br', \bl', \tr', \tl', \chi', \cdot')$ by
$(\lambda, \rho)$ as follows:
\begin{align}
& \qquad\rho(a\br v)=a\br'\rho(v), \qquad\qquad
a\tl v+\lambda(a\br v)=a\lambda(v)+a\tl'\rho(v),               \label{uc1}\\
& \qquad\rho(v\bl a)=\rho(v)\bl' a, \qquad\qquad
v\tr a+\lambda(v\bl a)=\lambda(v)a+\rho(v)\tr'a,               \label{uc2}\\
& \qquad\rho(v_{1}\cdot v_{2})=\rho(v_{1})\cdot'\rho(v_{2})
+\lambda(v_{1})\br'\rho(v_{2})+\rho(v_{1})\bl'\lambda(v_{2}),     \label{uc3}\\
& \chi(v_{1}, v_{2})+\lambda(v_{1}\cdot v_{2})
=\lambda(v_{1})\lambda(v_{2})+\rho(v_{1})\tr'\lambda(v_{2})
+\lambda(v_{1})\tl'\rho(v_{2})+\chi'(\rho(v_{1}), \rho(v_{2})),   \label{uc4}
\end{align}
for any $a\in A$, $v_{1}, v_{2}\in V$, then $\Omega(A, V)$ and  $\Omega'(A, V)$ are
called {\rm equivalent} and denoted by $\Omega(A, V)\approx\Omega'(A, V)$.
Moreover, if $\rho=\id$, $\Omega(A, V)$ and  $\Omega'(A, V)$ are called {\rm cohomologous}
denoted by $\Omega(A, V)\equiv\Omega'(A, V)$.
\end{defi}

It is easy to see that the relations ``$\approx$" and ``$\equiv$" are equivalent relations
on the set $\mathfrak{E}(A, V)$. Moreover, we can show that these relations corresponding
to the ``equivalent'' and ``cohomologous" relations defined in Definition \ref{Def:exte-equ}.

\begin{lem}\label{Lem:unif-coho}
Let $\Omega(A, V)=(\br, \bl, \tr, \tl, \chi, \cdot)$ and $\Omega'(A, V)=(\br', \bl',
\tr', \tl', \chi', \cdot')$ be two perm extending structures of $A$ by $V$,
$A\natural V$ and $A\natural' V$ be the corresponding unified products respectively.
Then, there is a bijection between the set of all morphisms of perm algebras
$\varphi: A\natural V\rightarrow A\natural' V$ which stabilizes $A$ and the set of
pairs $(\lambda, \rho)$, where $\lambda: V\rightarrow A$, $\rho: V\rightarrow V$ are
linear maps satisfying Eqs. (\ref{uc1})-(\ref{uc4}). More precisely,
For a pair $(\lambda, \rho)$, the corresponding homomorphism of perm algebras
$\varphi=\varphi_{(\lambda, \rho)}: A\natural V\rightarrow A\natural' V$ is given by
$\varphi(a, v)=(a+\lambda(v),\; \rho(v))$.

Moreover, The morphism $\varphi=\varphi_{(\lambda, \rho)}$ is an isomorphism
if and only if $\rho$ is a bijection, and $\varphi=\varphi_{(\lambda, \rho)}$
co-stabilizes $V$ if and only if $\rho=\id$.
\end{lem}

\begin{proof}
If $\varphi: A\natural V\rightarrow A\natural' V$ is a linear map and stabilizes $A$.
Then $\varphi(a, 0)=(a, 0)$ for any $a\in A$, and $\varphi(0, v)=(\lambda(v), \rho(v))$ with
some linear maps $\lambda: V\rightarrow A$ and $\rho: V\rightarrow V$. That is,
$\varphi(a, v)=(a+\lambda(v),\; \rho(v))$ and $\varphi$ is uniquely determined by
$(\lambda, \rho)$. Now we prove that $\varphi$ is a morphism of perm
algebras if and only if Eqs. (\ref{uc1})-(\ref{uc4}) hold. For any
$(a_{1}, v_{1}), (a_{2}, v_{2})\in A\oplus V$, note that
\begin{align*}
&\; \varphi\Big((a_{1}, v_{1})\ast(a_{2}, v_{2})\Big)\\
=&\; \varphi\Big(a_{1}a_{2}+v_{1}\tr a_{2}+a_{1}\tl v_{2}
+\chi(v_{1}, v_{2}),\ \ v_{1}\cdot v_{2}+a_{1}\br v_{2}+v_{1}\bl a_{2}\Big)\\
=&\; \Big(a_{1}a_{2}+v_{1}\tr a_{2}+a_{1}\tl v_{2}+\chi(v_{1}, v_{2})
+\lambda(v_{1}\cdot v_{2})+\lambda(a_{1}\br v_{2})+\lambda(v_{1}\bl a_{2}),\\[-2mm]
&\qquad\qquad\qquad\qquad\qquad\qquad\qquad\quad\rho(v_{1}\cdot v_{2})
+\rho(a_{1}\br v_{2})+\rho(v_{1}\bl a_{2})\Big),
\end{align*}
and
\begin{align*}
&\; \varphi(a_{1}, v_{1})\ast'\varphi(a_{2}, v_{2})\\
=&\; \Big(a_{1}+\lambda(v_{1}),\; \rho(v_{1})\Big)\ast'
\Big(a_{2}+\lambda(v_{2}),\; \rho(v_{2})\Big)\\
=&\; \Big((a_{1}+\lambda(v_{1}))(a_{2}+\lambda(v_{2}))+\rho(v_{1})\tr'(a_{2}+\lambda(v_{2}))
+(a_{1}+\lambda(v_{1}))\tl'\rho(v_{2})+\chi'(\rho(v_{1}), \rho(v_{2})),\\[-2mm]
&\qquad\qquad\qquad\qquad\qquad\quad(a_{1}+\lambda(v_{1}))\br'\rho(v_{2})
+\rho(v_{1})\bl'(a_{1}+\lambda(v_{1}))+\rho(v_{1})\cdot'\rho(v_{2})\Big),
\end{align*}
we get $\varphi((a_{1}, v_{1})\ast(a_{2}, v_{2}))=\varphi(a_{1}, v_{1})\ast'
\varphi(a_{2}, v_{2})$ if and only if Eqs. (\ref{uc1})-(\ref{uc4}) hold.

If $\rho$ is bijection, $\varphi$ is an isomorphism of perm algebras with
the inverse $\varphi^{-1}(a, v)=(a-\lambda(\rho^{-1}(v)),\; \rho^{-1}(v))$ for all
$a\in A$ and $v\in V$. Conversely, if $\varphi$ is an isomorphism, one can verify that
$\rho$ is a bijection. Finally, it is obvious that $\varphi_{\lambda, \rho}$ co-stabilizes
$V$ if and only if $\rho=\id_{V}$. The proof is completed.
\end{proof}

By the above discussion, we get $A\natural V\approx A\natural' V$
if and only if $\Omega(A, V)\approx\Omega'(A, V)$, and $A\natural V\equiv A\natural' V$
if and only if $\Omega(A, V)\equiv\Omega'(A, V)$. Thus, the answer for the extending
structures problem of perm algebras is given by the following theorem.

\begin{thm}\label{Thm:ext-prob}
Let $A$ be a perm algebra, $E$ a vector space that contains $A$
as a subspace and $V$ a complement of $A$ in $E$. Then,
\begin{itemize}
\item[$(i)$] Denote $H^{2}(A, V):=\mathfrak{E}(A, V)/\approx$.
     Then, the map
     $$
     H^{2}(A, V)\rightarrow \mathrm{Ext}(E, A),\qquad\qquad
     \overline{\Omega(A,V)}\mapsto \overline{A\natural V},
     $$
     is a bijection, where $\overline{\Omega(A, V)}$ and $\overline{A\natural V}$
     are the equivalence classes of $\Omega(A, V)$ and $A\natural V$ under
     $``\approx"$, respectively;
\item[$(ii)$] Denote $\mathcal{H}^{2}(A, V):=\mathfrak{E}(A, V)/\equiv$.
     Then, the map
     $$
     \mathcal{H}^{2}(A, V)\rightarrow \mathbf{Ext}(E, A),\qquad\qquad
     [\Omega(A, V)]\mapsto [A\natural V],
     $$
     is a bijection, where $[\Omega(A, V)]$ and $[A\natural V]$ are the equivalence
     classes of $\Omega(A, V)$ and $A\natural V$ under $``\equiv"$, respectively.
\end{itemize}
\end{thm}

\subsection{Factorization problem for perm algebras} \label{subsec:fact}

Let $\Omega(A, V)=(\br, \bl, \tr, \tl, \chi, \cdot)$ be an extending datum of
a perm algebra $A$ by a vector space $V$. If $\chi$ is trivial, i.e., $\chi(v_{1}, v_{2})=0$
for any $v_{1}, v_{2}\in V$, then $\Omega(A, V)=(\br, \bl, \tr, \tl, \chi, \cdot)$ be a perm
extending structure of $A$ by $V$ if and only if $(V, \cdot)$ is a perm algebra,
$(A, \tr, \tl)$ is a $V$-bimodule, and they satisfy Eqs. (\ref{ext1})-(\ref{ext9}).
In this case, the unified product $A\natural V$ denoted by $A\bowtie V$, is called
the {\it bicrossed product} of $A$ by $V$ and $(A, V, \br, \bl, \tr, \tl)$ is called
the {\it matched pair} of perm algebras $A$ and $V$. The bicrossed product of two perm
algebras is the construction responsible for the so-called {\it factorization problem}.

\vspace{2mm}
\noindent {\bf Factorization problem}: {\it Let $A$ and $V$ be two given perm algebras.
Describe and classify all perm algebras $E$ that factorize through $A$ and $V$, i.e.,
$E$ contains $A$ and $V$ as perm subalgebras such that $E=A+V$ and $A\cap V=\{0\}$.}

\vspace{2mm}
By the Theorem \ref{Thm:ext-unif}, we get the following corollary.

\begin{cor}\label{Cor:bicr}
A perm algebra $E$ factorizes through two given perm algebras $A$ and $V$ if and only if
there exists a matched pair of perm algebras $(A, V, \br, \bl, \tr, \tl)$ such that
$E\cong A\bowtie V$.
\end{cor}

\begin{proof}
If $E\cong A\bowtie V$, note that $A\cong A\bowtie\{0\}$ and $V \cong\{0\}\bowtie V$ are
perm subalgebras of $A \bowtie V $ and of course $A\bowtie V $ factorizes through
$A\times \{0\}$ and $\{0\}\times V$, we get $E$ factorizes through $A$ and $V$.
Conversely, if a perm algebra $E$ factorizes through two perm subalgebras $A$ and $V$.
Since $V$ is a subalgebra of $E$ and $A\cap V=\{0\}$, the cocycle $\chi: V\times V\rightarrow
A$ constructed in the proof of Theorem \ref{Thm:ext-unif} is a trivial map. Thus, the
unified product $A\natural V=A\bowtie V$ coincides with the bicrossed product of the perm
algebras $A$ and $V=\Ker(\pi_{A})$, and the map $\varphi: A\bowtie V\rightarrow E$,
$\varphi(b, v)=b+v$, becomes an isomorphism of perm subalgebras.
\end{proof}

Let $A$ be a perm subalgebra of $E$. A perm subalgebra $V$ of $E$ is called a {\it complement
of $A$ in $E$} (or $A$-complement of $E$) if $E=A+V$ and $A\cap V=\{0\}$. Based on Crocollary
\ref{Cor:bicr} we can restate the factorization problem as follows:
Let $A$ and $V$ be two given perm algebras. Describe the set of all matched
pairs $(A, V, \br, \bl, \tr, \tl)$ and classify up to an isomorphism all bicrossed
products $A\bowtie V$. For a perm subalgebra $A$ of $E$, denote
$\mathcal{F}(A, E)$ the set of the isomorphism classes of all $A$-complements in $E$.
Define the factorization index of $A$ in $E$ as $[E: A]=|\mathcal{F}(A, E)|$.

\begin{defi} \label{Def:unif-coho}
Let $(A, V, \br, \bl, \tr, \tl)$ be a matched pair of perm algebras.
A linear map $\psi: V\rightarrow A$ is called a {\rm deformation map} of the matched
pair $(A, V, \br, \bl, \tr, \tl)$ if for any $v_{1}, v_{2}\in V$,
\begin{align*}
&\; \psi(\mu_{V}(v_{1}, v_{2}))-\mu_{V}(\psi(v_{1}), \psi(v_{2}))\\
=&\; v_{1}\tr\psi(v_{2})+\psi(v_{1})\tl v_{2}-\psi(\psi(v_{1})\br v_{2})
-\psi(v_{1}\bl\psi(v_{2})).
\end{align*}
We denote the set of all deformation maps of the matched pair $(A, V, \br, \bl, \tr, \tl)$
by $\mathbf{DM}_{\br, \bl, \tr, \tl}(A, V)$.
\end{defi}

Now we study the classifying complements problem for perm algebras using
the concept of deformation map.

\begin{thm}\label{Thm: def-map}
Let $A$ be a perm subalgebra of $E$, $V$ be a given $A$-complement of $E$ with
the associated matched pair $(A, V, \br, \bl, \tr, \tl)$.
\begin{itemize}
\item[$(i)$] Let $\psi: V\rightarrow A$ be a deformation map of above matched pair. Then,
     the vector space $V$ with the new product:\\[-6mm]
     $$
     v_{1}\cdot_{\psi}v_{2}=\mu_{V}(v_{1}, v_{2})+\varphi(v_{1})\br v_{2}
     +v_{1}\bl\varphi(v_{1}),
     $$
     is a perm algebra, which is denoted by $V_{\psi}$. It is easy to see that
     $V_{\psi}$ is a $A$-complement of $E$, which is called the $\psi$-deformation of $V$.
\item[$(ii)$] $V'$ is a $A$-complement of $E$ if and only if $V'$ is isomorphic to
     $V_{\psi}$ for some deformation map $\psi: V\rightarrow A$ of the matched pair
     $(A, V, \br, \bl, \tr, \tl)$.
\end{itemize}
\end{thm}

\begin{proof}
$(i)$ Thank to Crocollary \ref{Cor:bicr}, we identity $E$ with $A\bowtie V$. For a
deformation map $\psi: V\rightarrow A$, we define $f_{\psi}: V\rightarrow E=A\bowtie V$
by $f_{\psi}(v)=(\psi(v), v)$ for any $v\in V$. Since $\psi$ is a deformation map,
for any $v_{1}, v_{2}\in V$, we have
\begin{align*}
&\; (\psi(v_{1}), v_{1})\ast(\psi(v_{2}), v_{2})\\
=&\; \Big(\psi(v_{1})\psi(v_{2})+v_{1}\tr\psi(v_{2})+\psi(v_{1})\tl v_{2},\ \
\mu_{V}(v_{1}, v_{2})+\psi(v_{1})\br v_{2}+v_{1}\bl\psi(v_{2})\Big)\\
=&\; \Big(\psi(\mu_{V}(v_{1}, v_{2})+\psi(v_{1})\br v_{2}+v_{1}\bl\psi(v_{2})),\ \
\mu_{V}(v_{1}, v_{2})+\psi(v_{1})\br v_{2}+v_{1}\bl\psi(v_{2})\Big).
\end{align*}
This means $f_{\psi}(v_{1})\ast f_{\psi}(v_{2})\in\Img(f_{\psi})$ for any $v_{1}, v_{2}\in V$,
and so that $\Img(f_{\psi})$ is a perm subalgebra of $E=A\bowtie V$. We identity $A$ with
$A\bowtie\{0\}$, and view it as a perm subalgebra of $E=A\bowtie V$. Then $A\cap\Img(f_{\psi})
=\{0\}$, and for any $(a, v)\in A\bowtie V$, $(a, v)=(a-\psi(v), 0)+(\psi(v), v)\in A+V$.
Hence, $\Img(f_{\psi})$ is a $A$-complement of $E=A\bowtie V$. Next, we show that
$\Img(f_{\psi})$ and $V_{\psi}$ are isomorphic as perm algebras. Denote $g_{\psi}:
V\rightarrow\Img(f_{\psi})$, $v\mapsto(\psi(v), v)$ for all $v\in V$. Then $g_{\psi}$ is a
linear bijection, and by the calculation as above, for any $v_{1}, v_{2}\in V$, we have
$$
g_{\psi}(v_{1}\cdot_{\psi}v_{2})=f_{\psi}(v_{1}\cdot_{\psi}v_{2})
=f_{\psi}(v_{1})\ast f_{\psi}(v_{2})=g_{\psi}(v_{1})\ast g_{\psi}(v_{2}).
$$
Thus, $\Img(f_{\psi})\cong V_{\psi}$ as perm algebras.

$(ii)$ Let $V'$ be an arbitrary of $A$-complement of $E$. Since $E=A\oplus V=A\oplus V'$
as vector spaces, there are four linear maps $\alpha: V\rightarrow A$,
$\beta: V\rightarrow V'$, $\epsilon: V'\rightarrow A$ and $\varepsilon: V'\rightarrow V$
such that
$$
v=\alpha(v)+\beta(v),\qquad\qquad v'=\epsilon(v')+\varepsilon(v'),
$$
for any $v\in V$ and $v'\in V'$. One can check that $\beta$ is a bijection.
We denote $\tilde{\beta}=\iota\circ\beta: V\rightarrow A\bowtie V$ the composition $\beta$ and
$\iota: V'\hookrightarrow E=A\bowtie V$. Then $\tilde{\beta}(v)=(-\alpha(v), v)$ for $v\in V$.
Moreover, since $V'=\Img(\beta)=\Img(\tilde{\beta})$ is a perm subalgebra
of $E$, we get for any $v_{1}, v_{2}\in V$, there exists a $v_{3}\in V$ such that
$(-\alpha(v_{1}), v_{1})\ast(-\alpha(v_{2}), v_{2})=(-\alpha(v_{3}), v_{3})$. This means
\begin{align*}
-\alpha(v_{3})&=\alpha(v_{1})\alpha(v_{2})-v_{1}\tr\alpha(v_{2})-\alpha(v_{1})\tl v_{2},\\
v_{3}&=\mu_{V}(v_{1}, v_{2})-\alpha(v_{1})\br v_{2}-v_{1}\bl\alpha(v_{2}).
\end{align*}
Applying $\alpha$ to the both sides of the first equation, and comparing with
the second equation, we can find that $-\alpha$ is a deformation map of the
matched pair $(A, V, \br, \bl, \tr, \tl)$. Finally, by direct calculation,
one can check that $\beta: V_{-\alpha}\rightarrow V'$ is an isomorphism of perm algebras.
\end{proof}

\begin{defi} \label{Def: def-equ}
Let $(A, V, \br, \bl, \tr, \tl)$ be a matched pair of perm algebras.
For two deformation maps $\psi$, $\phi: V\rightarrow A$, if there exists a linear
automorphism $\rho: V\rightarrow V$ such that for any $v_{1}, v_{2}\in V$,
\begin{align*}
&\; \rho(\mu_{V}(v_{1}, v_{2}))-\mu_{V}(\rho(v_{1}), \rho(v_{2}))\\
=&\; \phi(\rho(v_{1}))\br\rho(v_{2})+\rho(v_{1})\bl\phi(\rho(v_{2}))
-\rho(\psi(v_{1})\br v_{2})-\rho(v_{1}\bl\phi(v_{2})),
\end{align*}
then $\psi$ and $\phi$ are called {\rm equivalent}, and denoted by $\psi\sim\phi$.
\end{defi}

\begin{thm}\label{Thm: def-map}
Let $A$ be a perm subalgebra of $E$, $V$ be a $A$-complement of $E$ and $(A, V, \br, \bl, \tr,
\tl)$ be the associated matched pair. Then, $\sim$ is an equivalence relation on the set
$\mathbf{DM}(A, V, \br, \bl$, $\tr, \tl)$. Denote $\mathcal{H}^{2}_{\br, \bl, \tr, \tl}(A, V)
=\mathbf{DM}_{\br, \bl, \tr, \tl}(A, V)/\sim$. The map
$$
\mathcal{H}^{2}_{\br, \bl, \tr, \tl}(A, V)\rightarrow\mathcal{F}(A, E), \qquad\qquad
[\psi]\mapsto \bar{V}_{\psi},
$$
is a bijection, where $[\psi]$ is the equivalence class of $\psi$, $\bar{V}_{\psi}$
is the isomorphism class of $V_{\psi}$.
In particular, $[E: A]=|\mathcal{H}^{2}_{\br, \bl, \tr, \tl}(A, V)|$.
\end{thm}

\begin{proof}
By direct calculation, we botain $\psi$ and $\phi$ are called equivalent if and only if
there exists a linear bijection $\rho: V\rightarrow V$ such that for any $v_{1}, v_{2}\in V$,
\begin{align*}
&\; \rho(v_{1}\cdot_{V} v_{2})+\rho(\psi(v_{1})\br v_{2})+\rho(v_{1}\bl\phi(v_{2}))\\
=&\; \rho(v_{1})\cdot_{V}\rho(v_{2})+\phi(\rho(v_{1}))\br\rho(v_{2})
+\rho(v_{1})\bl\phi(\rho(v_{2})),
\end{align*}
if and only if $\rho: V_{\psi}\rightarrow V_{\phi}$ is an isomorphism
of perm algebras. We get this theorem.
\end{proof}

\subsection{Flag extending structures} \label{subsec: flag}
In this subsection, we apply our main theorem to a flag extending structure of
perm algebras as a special case, and provide a way of answering the problem
for a special case.

\begin{defi} \label{Def:flag}
Let $(A, R)$ be a perm algebra and $E$ be a vector space containing $A$ as a subspace.
A perm algebra structure $\diamond$ on $E$ such that $A$ is a perm subalgebra
is called a {\rm flag extending structure} of $A$ to $(E, \diamond)$ if there is
a finite chain of perm subalgebras of $(E, \diamond)$:
$$
A=E_{0}\subset E_{1}\subset\cdots\subset E_{m}=E
$$
such that $E_{i}$ has codimension 1 in $E_{i+1}$, for all $i=0, 1,\cdots, m-1$.
\end{defi}

In this case, if the complement of $A$ in $E$ is finite-dimensional,
we can obtain all flag extending structure of $A$ to $E$.
In fact, all flag extending structures of $A$ to $E$ can be completely
described by a recursive process. First, classify all unified products of $A\natural V_{1}$,
where $V_{1}$ is a 1-dimensional vector space. Second, replacing $A$ by
$A\natural V_{1}$, characterize all unified products of $(A\natural V_{1})\natural V_{2}$
where $V_{2}$ is a 1-dimensional vector space. After finite steps, we can describe and
classify all flag extending structures of $A$ to $E$.
Therefore, in this section, we mainly study the extending structure of $A$
by a 1-dimensional vector space $V$. For this special extending structure,
we introduce the notion of flag datum.

\begin{defi}\label{Def:fla-dat}
Let $A$ be a perm algebra. A {\rm flag datum} of $A$ is a 6-tuple $(h, g, D, T, \tilde{a},
\tilde{k})$ where $\tilde{a}\in A$, $\tilde{k}\in\bk$, $g, h: A\rightarrow\bk$,
$D, T: A\rightarrow A$ are linear maps satisfying for $a, a_{1}, a_{2}\in A$:
\begin{align}
h(a_{1}a_{2})&=h(a_{1})h(a_{2})=h(a_{1})g(a_{2}), \qquad\quad h(T(a))=0,   \label{es1}\\
g(a_{1}a_{2})&=g(a_{2}a_{1})=g(a_{1})g(a_{2}),                            \label{es2}\\
D(a_{1}a_{2})&=D(a_{2}a_{1})=D(a_{1})a_{2}+g(a_{1})D(a_{2}),             \label{es3}\\
T(a_{1}a_{2})&=a_{1}T(a_{2})+T(a_{1})h(a_{2})
=a_{1}D(a_{2})+T(a_{1})g(a_{2})=T(a_{1})a_{2}+h(a_{1})D(a_{2}),          \label{es4}\\
\tilde{k}g(a)&=\tilde{k}g(a)+g(D(a))=\tilde{k}h(a)+g(T(a))
=\tilde{k}g(a)+h(D(a)),                                                 \label{es5}\\
a\tilde{a}+\tilde{k}T(a)&=T(T(a))+h(a)\tilde{a},                        \label{es6}\\
a\tilde{a}+\tilde{k}D(a)&=D(D(a))+g(a)\tilde{a}=D(T(a))+h(a)\tilde{a}
=T(D(a))+g(a)\tilde{a},                                                  \label{es7}\\
g(\tilde{a})&=h(\tilde{a})\qquad\quad T(\tilde{a})=D(\tilde{a}).        \label{es8}
\end{align}
Denote by $\mathcal{F\!\!L}(A)$ the set of all flag datums of $A$.
\end{defi}

\begin{pro}\label{Pro:e-f}
Let $(A, R)$ be a perm algebra and $V$ be a vector space of dimension 1 with a basis $\{x\}$.
Then, there exists a bijection between the set $\mathfrak{E}(A, V)$ of all perm
extending structures of $A$ by $V$ and $\mathcal{F\!\!L}(A)$.

Through the above bijection, the perm extending structure $\Omega(A, V)=(\br, \bl,
\tr, \tl, \chi, \cdot)$ corresponding to $(h, g, D, T, \tilde{a}, \tilde{k})\in
\mathcal{F\!\!L}(A)$ is given as follows:
\begin{align*}
& a\br x=h(a)x,\qquad\, x\bl a=g(a)x,\qquad\;\, x\tr a=D(a),\\
& \; a\tl x=T(a),\qquad\quad x\cdot x=\tilde{k}x,\qquad\quad \chi(x, x)=\tilde{a},
\end{align*}
for all $a\in A$.
\end{pro}

\begin{proof}
Given a perm extending structure $\Omega(A, V)=(\br, \bl, \tr, \tl, \chi, \cdot)$.
Since $V=\bk\{x\}$, we can set
\begin{align*}
& a\br x=h(a)x,\qquad\, x\bl a=g(a)x,\qquad\;\, x\tr a=D(a),\\
& \; a\tl x=T(a),\qquad\quad x\cdot x=\tilde{k}x,\qquad\quad \chi(x, x)=\tilde{a},
\end{align*}
where $\tilde{a}\in A$, $\tilde{k}\in\bk$, $g, h: A\rightarrow\bk$ and
$D, T: A\rightarrow A$ are linear maps. By a straightforward computation,
we can obtain that the conditions $(\ref{ext1})$-$(\ref{ext10})$ in Proposition
\ref{Pro:exten} are equivalent to the fact that $g: A\rightarrow\bk$ is a
morphism of algebras and $(\ref{es1})$-$(\ref{es8})$ hold.
\end{proof}

Given a flag datum $\Gamma=(h, g, D, T, \tilde{a}, \tilde{k})$ of
$A$, in the proposition, there is a unique perm extending structure
corresponding to it. We denote by this extending structure by $\Omega_{\Gamma}(A, V)$.
Then the one-to-one correspondence given in Proposition \ref{Pro:e-f} is just
$$
\Phi:\; \mathcal{F\!\!L}(A)\rightarrow\mathfrak{E}(A, V),\qquad
\Gamma=(h, g, D, T, \tilde{a}, \tilde{k})\mapsto \Omega_{\Gamma}(A, V).
$$
Next, we consider an equivalent relation on $\mathcal{F\!\!L}(A)$.

\begin{defi}\label{Def:fla-equ}
Let $(h, g, D, T, \tilde{a}, \tilde{k})$ and $(h', g', D', T', \tilde{a}', \tilde{k}')$
be two flag datums of a perm algebra $A$. They are called {\rm equivalent}, denoted by
$(h, g, D, T, \tilde{a}, \tilde{k})\approx(h', g', D', T', \tilde{a}', \tilde{k}')$, if
$h=h'$, $g=g'$, and there exists an element $\vec{a}\in A$ and a
nonzero element $l\in\bk$ such that for any $a\in A$,
\begin{align}
D(a)&=\vec{a}a+lD'(a)-g(a)\vec{a},                               \label{fq1}\\
T(a)&=a\vec{a}+lT'(a)-h(a)\vec{a},                               \label{fq2}\\
\tilde{a}&=\vec{a}^{2}+lD'(\vec{a})+lT'(\vec{a})
+l^{2}\tilde{a}'-\tilde{k}\vec{a},                               \label{fq3}\\
\tilde{k}&=l\tilde{k}'+h'(\vec{a})+g'(\vec{a}),                  \label{fq4}
\end{align}
\end{defi}

Now, we can use the flag datums to classify some special extensions of perm algebras.

\begin{thm}\label{Thm:fla-bij}
Let $A$ be a perm algebra, $E$ be a vector space that contains $A$ as a subspace and
$V=\bk\{x\}$ be a complement of $A$ in $E$. Then the relation $``\approx"$ is an
equivalent relation of the set $\mathcal{F\!\!L}(A)$, and there is a bijection between
$\mathrm{Ext}(E, A)$ and $\mathcal{F\!\!L}(A)/\approx$.
\end{thm}

\begin{proof}
First, by a straightforward verification shows that $\approx$ is an
equivalent relation of the set $\mathcal{F\!\!L}(A)$. Second, thank to Theorem
\ref{Thm:ext-prob}, we only need show that $\Phi$ induces a one-to-one correspondence
between $H^{2}(A, V)$ and $\mathcal{F\!\!L}(A)/\approx$. For any two flag datums
$\Gamma=(h, g, D, T, \tilde{a}, \tilde{k})$ and $\Gamma'(h', g', D', T', \tilde{a}',
\tilde{k}')$ of $A$, we denote the corresponding extending structures by
$\Omega_{\Gamma}(A, V)=(\br$, $\bl, \tr, \tl, \chi, \cdot)$ and $\Omega_{\Gamma'}(A, V)
=(\br', \bl', \tr', \tl', \chi', \cdot')$ respectively. By a straightforward computation,
we can obtain that there is a linear map $\lambda: V\rightarrow A$ and a linear
bijiction $\rho: V\rightarrow V$, such that the Eqs. $(\ref{uc1})$-$(\ref{uc4})$ hold,
is equivalent to, there exist an $\vec{a}\in A$ and a nonzero element $l\in\bk$ such
that Eqs. $(\ref{fq1})$-$(\ref{fq4})$ hold. It is necessary to point out that
Eqs. $(\ref{uc1})$ and $(\ref{uc2})$ imply $h=h'$ and $g=g'$.
\end{proof}

After provid a theoretical answer to the extending structures problem,
we are left to give the extension algebra $E$ by calculation for a given perm algebra 
$A$ and vector space $V$. In general, it is quite difficult. But when $\dim_{\bk}V=1$, 
we can obtain all the extension algebras $E$ by compute $\mathcal{F\!\!L}(A)$.

\begin{ex}\label{ex:fla-bij}
Let $\bk=\bc$ be the field of complex numbers. We consider a 2-dimensional
perm $\bk$-algebra $A$ with a basis $\{e_{1}, e_{2}\}$, the non-zero products
$e_{1}e_{1}=e_{1}$ and $e_{2}e_{1}=e_{2}$. By direct computations, we can obtain
that the flag datum $(h, g, D, T, \tilde{a}, \tilde{k})$ of $A$ is given by as follows:
\begin{align*}
g(e_{1})&=k,\qquad h(e_{1})=l,\qquad\; D(e_{1})=se_{2},\qquad T(e_{1})=te_{2},\\
g(e_{2})&=0,\qquad h(e_{2})=0,\qquad D(e_{2})=0,\qquad\quad T(e_{2})=0,\qquad \tilde{a}=qe_{2}
\end{align*}
where $k, l, s, t, q\in\bk$, and $k=0$ or $1$, $l=0$ or $1$, $kl=l$, $ks=ls=(k-1)y
=\tilde{k}(k-l)=0$, $\tilde{k}t=lq$, $\tilde{k}s=kq=lq$.
Therefore, we can classify the flag datum by the following case.

{\bf Case I.} $k=l=1$. In this case, we have $s=q=0$ and $\tilde{k}t=0$.
$(i)$ if $\tilde{k}=0$, the perm algebra structure on $\bc\{e_{1}, e_{2}, x\}$ induced
by the flag datum is given by $e_{1}e_{1}=e_{1}$, $e_{1}x=te_{2}+x$, $e_{2}e_{1}=e_{2}$,
$xe_{1}=x$, for $t\in\bc$, where we only give the nonzero produces.
$(ii)$ if $t=0$, the perm algebra structure on $\bc\{e_{1}, e_{2}, x\}$ induced by the
flag datum is given by $e_{1}e_{1}=e_{1}$, $e_{1}x=x$, $e_{2}e_{1}=e_{2}$,
$xe_{1}=x$, $xx=\tilde{k}x$, for $\tilde{k}\in\bc$. By Theorem \ref{Thm:fla-bij},
we get each of these perm algebras is isomorphic to one of the following perm algebra
on $\bc\{e_{1}, e_{2}, x\}$:
\begin{align*}
&A_{1}:\qquad e_{1}e_{1}=e_{1},\quad e_{1}x=xe_{1}=x,\quad e_{2}e_{1}=e_{2};\\
&A_{2}:\qquad e_{1}e_{1}=e_{1},\quad e_{1}x=xe_{1}=x,\quad e_{2}e_{1}=e_{2},\quad xx=x.
\end{align*}

{\bf Case II.} $k=1$ and $l=0$. In this case, we have $s=q=\tilde{k}=0$, and
the perm algebra structure on $\bc\{e_{1}, e_{2}, x\}$ induced by the flag datum
is given by $e_{1}e_{1}=e_{1}$, $e_{1}x=te_{2}$, $e_{2}e_{1}=e_{2}$,
$xe_{1}=x$, for $t\in\bc$. By Theorem \ref{Thm:fla-bij},
we get each of these perm algebras is isomorphic to one of the following perm algebra
on $\bc\{e_{1}, e_{2}, x\}$:
\begin{align*}
&A_{3}:\qquad e_{1}e_{1}=e_{1},\quad e_{1}x=e_{2},\quad e_{2}e_{1}=e_{2},\quad xe_{1}=x;\\
&A_{4}:\qquad e_{1}e_{1}=e_{1},\quad e_{2}e_{1}=e_{2},\quad xe_{1}=x.
\end{align*}

{\bf Case III.} $k=l=0$. In this case, we have $t=0$ and $\tilde{k}s=0$.
$(i)$ if $\tilde{k}=0$, the perm algebra structure on $\bc\{e_{1}, e_{2}, x\}$ induced
by the flag datum is given by $e_{1}e_{1}=e_{1}$, $e_{1}x=te_{2}$, $e_{2}e_{1}=e_{2}$,
$xe_{1}=x$, $xx=qe_{2}$, for $t, q\in\bc$. By Theorem \ref{Thm:fla-bij},
we get each of these perm algebras is isomorphic to one of the following perm algebra
on $\bc\{e_{1}, e_{2}, x\}$:
\begin{align*}
&A_{5}:\qquad e_{1}e_{1}=e_{1},\quad e_{2}e_{1}=e_{2}, \quad xe_{1}=x,\quad xx=e_{2};\\
&A_{6}:\qquad e_{1}e_{1}=e_{1},\quad e_{2}e_{1}=e_{2}, \quad xe_{1}=x.
\end{align*}
$(ii)$ if $s=0$, the perm algebra structure on $\bc\{e_{1}, e_{2}, x\}$ induced by the
flag datum is given by $e_{1}e_{1}=e_{1}$, $e_{2}e_{1}=e_{2}$, $xx=qe_{2}+\tilde{k}x$,
for $q, \tilde{k}\in\bc$. We denote these perm algebras by $A_{q, \tilde{k}}$,
$q, \tilde{k}\in\bc$. By Theorem \ref{Thm:fla-bij}, we get $A_{q, \tilde{k}}\cong
A_{q', \tilde{k}'}$ if and only if there exist a nonzero $p\in\bc$, such that
$\tilde{k}=p\tilde{k}'$ and $q=p^{2}q'$. This means there are infinite many perm algebras
that are not isomorphic to each other.

Therefore, we get all the perm algebra structure on $\bc\{e_{1}, e_{2}, x\}$, such that
$A$ is a perm sualgebra of it.
\end{ex}

\section{Non-abelian extensions and automorphisms of Perm algebras} \label{sec: non-abel}

In this section, we apply the our results to the non-abelian extensions of perm
algebras and provide the the fundamental sequence of Wells in the context of perm algebras.
Let $\Omega(A, V)=(\br, \bl, \tr, \tl, \chi, \cdot)$ be a perm extending structure.
If the actions $``\br"$ and $``\bl"$ are trivial, then the corresponding unified product
$A\natural V$ is also called a {\it crossed product} of perm algebra $A$ by $V$.
This crossed product is closely related to the non-abelian extension of perm algebras.

\subsection{Non-abelian extension and crossed product} \label{subsec: abl-prod}
Here, we define a non-abelian cohomology group of perm algebras,
and show that the non-abelian extensions can be classified by the non-abelian cohomology.
First, we give the definition of non-abelian extensions of perm algebras.

\begin{defi}\label{def:nonabel}
Let $A$ and $B$ be two perm algebras. A {\rm non-abelian extension} $\mathcal{E}$ of
$B$ by $A$ is a short exact sequence in the category of perm algebras as follows:
$$
\mathcal{E}:\qquad \xymatrix{0\ar[r]& A\ar[r]^{\iota}
& E\ar[r]^{\pi}& B\ar[r]&0,}
$$
where $E$ is called a {\rm non-abelian extension perm algebra} of $B$ by $A$.
In particular, if the multiplication of $A$ is trivial, i.e., $a_{1}a_{2}=0$ for
any $a_{1}, a_{2}\in A$, then $\mathcal{E}$ is called an {\rm abelian extension
perm algebra} of $B$ by $A$.

Let $\mathcal{E}$ and $\mathcal{E}'$ two extensions of $B$ by $A$. They are called
{\rm equivalent} if there exists a morphism of perm algebras $\varphi: E\rightarrow E'$
such that the following diagram commutes
$$
\xymatrixrowsep{0.5cm}
\xymatrix{\mathcal{E}:\qquad 0\ar[r]& A\ar[r]^{\iota}
\ar@{=}[d]& E\ar[r]^{\pi}\ar[d]^{\varphi} & B\ar[r]\ar@{=}[d]&0\;\\
\mathcal{E}':\qquad 0\ar[r]& A\ar[r]^{\iota'} & E'\ar[r]^{\pi'}& B\ar[r]&0.}
$$
We denote by $\mathcal{E}xt(B, A)$ the set of all equivalence classes of non-abelian
extension of $B$ by $A$, and by $[\mathcal{E}]$ the equivalence class of extension
$\mathcal{E}$.
\end{defi}

Let $\mathcal{E}:\xymatrix@C=0.5cm{0\ar[r] & A\ar[r]^{\iota}& E \ar[r]^{\pi}& B\ar[r]&0}$
be a non-abelian extension of perm algebras. We often identify $A$ as a perm
subalgebra of $E$. Since $\mathcal{E}$ is split in the category of vector spaces,
there exists a linear map $\varsigma: B\rightarrow E$ such that $\pi\circ\varsigma
=\id_{B}$. We call $\varsigma$ a section of $\mathcal{E}$. Then there is a $B$ bimodule
structure on $A$ by the multiplication of $E$, i.e., $b\tr a=\varsigma(b)\diamond a$ and
$a\tr b=a\diamond\varsigma(b)$ for $a\in A$, $b\in B$, and this bimodule does not depend
on the choice of $\varsigma$.

\begin{defi}\label{def:coho}
Let $A$ and $B$ be two perm algebras. A {\rm non-abelian 2-cocycle} on $B$ with values in
$A$ is a triple $(\tr, \tl, \chi)$ of bilinear maps $\tr: B\times A\rightarrow A$,
$\tl: A\times B\rightarrow A$ and $\chi: B\times B\rightarrow A$,
satisfying the following properties:
\begin{align}
(a_{1}a_{2})\tl b-(a_{2}a_{1})\tl b& =a_{1}(a_{2}\tl b)-a_{2}(a_{1}\tl b), \label{nab1}\\
b\tr(a_{1}a_{2})-a_{1}(b\tr a_{2})& =(b\tr a_{1}-a_{1}\tl b)a_{2},         \label{nab2}\\
b_{1}\tr (a\tl b_{2})-(b_{1}\tr a)\tl b_{2}& =
a\tl(b_{1}b_{2})-(a\tl b_{1})\tl b_{2}+a\chi(b_{1}, b_{2}),                \label{nab3}\\
(b_{1}b_{2}-b_{2}b_{1})\tr a& =b_{1}\tr(b_{2}\tr a)
-b_{2}\tr(b_{1}\tr a)-\chi(b_{1}, b_{2})a+\chi(b_{2}, b_{1})a,             \label{nab4}\\
\chi(b_{1}b_{2},\; b_{3})-\chi(b_{1},\; b_{2}b_{3})&
-\chi(b_{2}b_{1},\; b_{3})+\chi(b_{2},\; b_{1}b_{3})                    \label{nab5}\\[-1mm]
=b_{1}\tr\chi(b_{2}, b_{3})-b_{2}\tr\chi(b_{1}, &\, b_{3})
-\chi(b_{1}, b_{2})\tl b_{3}+\chi(b_{2}, b_{1})\tl b_{3},                   \nonumber
\end{align}
for any $a, a_{1}, a_{2}\in A$ and $b, b_{1}, b_{2}, b_{3}\in B$.
One denotes by $Z^{2}_{nab}(B, A)$ the set of theses cocycles.

Moreover, $(\tr, \tl, \chi)$ and $(\tr', \tl', \chi')$ are said to be
{\rm equivalent} if there exists a linear map $\lambda: B\rightarrow A$ satisfying:
\begin{align}
& \qquad\quad a\tl b-a\tl'b=a\lambda(b), \qquad\qquad
b\tr a-b\tr'a=\lambda(b)a,                                                    \label{eq1}\\
& \chi(b_{1}, b_{2})-\chi'(b_{1}, b_{2})=\lambda(b_{1})\lambda(b_{2})
+b_{1}\tr'\lambda(b_{2})+\lambda(b_{1})\tl'b_{2}-\lambda(b_{1}b_{2}),         \label{eq2}
\end{align}
for any $a\in A$ and $b, b_{1}, b_{2}\in B$. The {\rm non-abelian cohomology}
$H^{2}_{nab}(B, A)$ is the quotient of $Z^{2}_{nab}(B, A)$ by this equivalence relation.
\end{defi}

It is worth noting that $H^{2}_{nab}(B, A)$ is not a group in general.

\begin{pro} \label{pro: dir-sum}
Let $A$ and $B$ be two perm algebras, $\chi: B\times B\rightarrow A$,
$\tr: B\times A\rightarrow A$ and $\tl: A\times B\rightarrow A$ be bilinear maps.
Then the triple $(\tr, \tl, \chi)$ is a non-abelian 2-cocycle on $B$ with values
in $A$ if and only if the vector space $A\oplus B$ under the multiplication
$$
(a_{1}, b_{1})(a_{2}, b_{2})=\Big(a_{1}a_{2}+b_{1}\tr a_{2}+a_{1}\tl b_{2}
+\chi(b_{1}, b_{2}),\; b_{1}b_{2}\Big),
$$
is a perm algebra. In this case, we denote this perm algebra by $A\rtimes_{(\tr, \tl, \chi)}B$.
\end{pro}

Let $(\tr, \tl, \chi)$ be a non-abelian 2-cocycle on perm algebra $B$ with values in perm
algebra $A$. Then we get a non-abelian extension of $B$ by $A$ as follows:
$$
\mathcal{E}^{(\tr, \tl, \chi)}:\qquad \xymatrix{0\ar[r]& A\ar[r]^{\iota_{A}\qquad}
& A\rtimes_{(\tr, \tl, \chi)}B\ar[r]^{\qquad\pi_{B}}& B\ar[r]&0,}\qquad\qquad
$$
where $\iota_{A}$ and $\pi_{B}$ are canonical injection and projection.

Let $\Omega(A, V)=(\br, \bl, \tr, \tl, \chi, \cdot)$ be a perm extending structure.
If the actions $``\br"$ and $``\bl"$ are trivial, Eq. (\ref{ext10}) in
Proposition \ref{Pro:exten} means $(V, \cdot)$ is a perm algebra, Eqs.
(\ref{ext3}), (\ref{ext4}), (\ref{ext7})-(\ref{ext9}) hold is equivalent to
Eqs. (\ref{nab1})-(\ref{nab5}) in Definition \ref{def:coho} hold,
and Eqs. (\ref{eq1})-(\ref{eq2}) hold is equivalent to Eqs.
(\ref{uc1})-(\ref{uc4}) in Definition \ref{Def:unif-coho} hold for $\rho=\id$.
Now let $A$ and $B$ be two perm algebras. From the conclusions of subsection
\ref{subsec:unif}, we have
\begin{itemize}
\item[$(i)$] the perm algebra $A\rtimes_{(\tr, \tl, \chi)}B$
     just the crossed product of perm algebra of $A$ by $B$;
\item[$(ii)$] for each non-abelian extension algebra $E$ of $B$ by $A$, there exists a
     non-abelian 2-cocycle $(\tr, \tl, \chi)$ such that $E\cong A\rtimes_{(\tr, \tl, \chi)}B$
     as perm algebras, and for each non-abelian extension $\mathcal{E}$ of $B$ by $A$,
     there exists a non-abelian 2-cocycle $(\tr, \tl, \chi)$ such that $\mathcal{E}$
     and $\mathcal{E}^{(\tr, \tl, \chi)}$ are equivalent;
\item[$(iii)$] two non-abelian extensions $\mathcal{E}^{(\tr, \tl, \chi)}$ and
     $\mathcal{E}^{(\tr', \tl', \chi')}$ are equivalent if and only if non-abelian
     2-cocycles $(\tr, \tl, \chi)$ and $(\tr', \tl', \chi')$ are equivalent.
\end{itemize}

Thus, by Theorem \ref{Thm:ext-prob}, we have

\begin{thm} \label{thm:nonabel}
Let $A$ and $B$ be two perm algebras. There is a one-to-one correspondence
$$
\Psi:\quad H^{2}_{nab}(B, A)\rightarrow \mathcal{E}xt(B, A),
\qquad\qquad [(\tr, \tl, \chi)]\mapsto [\mathcal{E}^{(\tr, \tl, \chi)}],
$$
where $[(\tr, \tl, \chi)]$ is the the equivalent class of $(\tr, \tl, \chi)$.
In other words, $H^{2}_{nab}(B, A)$ classifies non-abelian extensions of $B$ by $A$.
\end{thm}

Indeed, for any non-abelian extension $\mathcal{E}:\xymatrix@C=0.5cm{0\ar[r]
& A \ar[r]^{\iota}& E \ar[r]^{\pi}& B\ar[r]&0}$ of perm algebras with a section $\varsigma$,
we can define bilinear maps $\tr^{\varsigma}: B\times A\rightarrow A$, $\tl^{\varsigma}:
A\times B\rightarrow A$ and $\chi^{\varsigma}: B\times B\rightarrow A$ by
\begin{align*}
b\tr^{\varsigma}a& =\varsigma(b)\diamond a, \qquad\qquad\qquad
a\tl^{\varsigma}b=a\diamond\varsigma(b),\\
&\chi^{\varsigma}(b_{1}, b_{2})=\varsigma(b_{1})\diamond\varsigma(b_{2})-\varsigma(b_{1}b_{2}),
\end{align*}
for any $a\in A$ and $b, b_{1}, b_{2}\in B$. Then one can check that
$(\tr^{\varsigma}, \tl^{\varsigma}, \chi^{\varsigma})$ is a non-abelian 2-cocycle,
which does not depend on the choice of $\varsigma$, and extensions $\mathcal{E}$ and
$\mathcal{E}^{(\tr^{\varsigma}, \tl^{\varsigma}, \chi^{\varsigma})}$ are equivalent.
This gives a inverse map of $\Psi$.

\subsection{The inducibility problem of automorphisms}\label{subsec: ind-aut}
In this subsection, we study the inducibility of a pair of automorphisms related to
a non-abelian extension of perm algebras, and give the fundamental sequence of Wells
in the context of perm algebras.

Let $A$ and $B$ be two perm algebras,
$$
\mathcal{E}:\quad \xymatrix{0\ar[r]& A\ar[r]^{\iota}& E\ar[r]^{\pi}& B\ar[r]&0}
$$
be a non-abelian extension of $B$ by $A$. We always regard $A$ as a subalgebra of $E$.
Denote
$$
\Aut_{A}(E):=\{\alpha\in\Aut(E)\mid\alpha(A)=A\}.
$$
Then, $\alpha$ preserves also $B$ and $\alpha|_{A}\in\Aut(A)$ if $\alpha\in\Aut_{A}(E)$.
Let $\varsigma: B\rightarrow E$ be a section of $\mathcal{E}$. For any
$\alpha\in\Aut_{A}(E)$, we define a linear map
$$
\bar{\alpha}:=\pi\circ\alpha\circ\varsigma:\quad  B\longrightarrow B.
$$
One can check that $\bar{\alpha}$ does not depend on the choice of the
section $\varsigma$. Since $\pi$ is the projection and $\alpha$ preserves $B$,
we get $\bar{\alpha}$ is a bijection. For any $b_{1}, b_{2}\in B$, we have
\begin{align*}
\bar{\alpha}(b_{1}b_{2})&=\pi(\alpha(\varsigma(b_{1}b_{2})))
=\pi(\alpha(\varsigma(b_{1})\varsigma(b_{2})-\chi(b_{1}, b_{2})))\\
&=\pi(\alpha(\varsigma(b_{1})\varsigma(b_{2})))
=\bar{\alpha}(b_{1})\bar{\alpha}(b_{2}).
\end{align*}
Thus $\bar{\alpha}\in\Aut(B)$. Note that $\overline{\beta\circ\alpha}
=\pi\circ\beta\circ\alpha\circ\varsigma=\pi\circ\beta\circ\varsigma\circ\pi\circ
\alpha\circ\varsigma=\bar{\beta}\circ\bar{\alpha}$, for any $\alpha, \beta\in\Aut_{A}(E)$,
we obtain a group homomorphism
$$
\kappa:\; \Aut_{A}(E)\rightarrow\Aut(A)\times\Aut(B),\qquad\qquad
\alpha\mapsto(\alpha|_{A},\; \bar{\alpha}).
$$

\begin{defi}\label{def:induc}
Let $\mathcal{E}: \xymatrix@C=0.5cm{0\ar[r]& A\ar[r]^{\iota}& E\ar[r]^{\pi}& B\ar[r]&0}$
be a non-abelian extension of perm algebra $B$ by perm algebra $A$. A pair
$(\beta, \gamma)\in\Aut(A)\times\Aut(B)$ of automorphisms is said to be
{\rm inducible} if there exists a map $\alpha\in\Aut_{A}(E)$ such that
$(\beta, \gamma)=(\alpha|_{A}, \bar{\alpha})$.
\end{defi}

Let $A$ and $B$ be two perm algebras. For any a pair $(\beta, \gamma)
\in\Aut(A)\times\Aut(B)$ and a non-abelian 2-cocycle $(\tr, \tl, \chi)$ on $B$ with
values in $A$, we define a new quadruple $(\tr^{\beta, \gamma}, \tl^{\beta, \gamma},
\chi^{\beta, \gamma})$ by
\begin{align*}
\tr^{\beta, \gamma}(b, a)&=\beta\Big(\gamma^{-1}(b)\tr \beta^{-1}(a)\Big),\\
\tl^{\beta, \gamma}(a, b)&=\beta\Big(\beta^{-1}(a)\tl \gamma^{-1}(b)\Big),\\
\chi^{\beta, \gamma}(b_{1}, b_{2})&=\beta(\chi(\gamma^{-1}(b_{1}),\; \gamma^{-1}(b_{2}))),
\end{align*}
for $a\in A$ and $b, b_{1}, b_{2}\in B$. Then one can check that
$(\tr^{\beta, \gamma}, \tl^{\beta, \gamma}, \chi^{\beta, \gamma})$ is also a non-abelian
2-cocycle on $B$ with values in $A$.

\begin{thm}\label{thm: induc}
Let $\mathcal{E}: \xymatrix@C=0.5cm{0\ar[r]& A\ar[r]^{\iota}& E\ar[r]^{\pi}& B\ar[r]&0}$
be a non-abelian extension of perm algebra $B$ by perm algebra $A$, $\varsigma$ be a
section of $\mathcal{E}$. With the above notations, a pair $(\beta, \gamma)\in\Aut(B)
\times\Aut(A)$ is inducible if and only if the non-abelian $2$-cocycles $(\tr, \tl, \chi)$
and $(\tr^{\beta, \gamma}, \tl^{\beta, \gamma}, \chi^{\beta, \gamma})$ are equivalent,
where $(\tr, \tl, \chi)$ is the non-abelian $2$-cocycle
corresponding to $\mathcal{E}$ with section $\varsigma$.
\end{thm}

\begin{proof}
If the pair $(\beta, \gamma)$ is inducible, there exists an element $\alpha\in\Aut_{A}(E)$
such that $\alpha|_{A}=\beta$ and $\bar{\alpha}=\pi\circ\alpha\circ\varsigma=\gamma$.
Note that for any $b\in B$, $\pi\Big(\varsigma(b)-\alpha\circ\varsigma\circ
\gamma^{-1}(b)\Big)=b-b=0$, we get $\varsigma(b)-\alpha\circ\varsigma\circ\gamma^{-1}(b)
\in\Ker(\pi)\cong A$. Thus we can define a linear map $\vartheta: B\rightarrow A$ by
$\vartheta(b)=\varsigma(b)-\alpha\circ\varsigma\circ\gamma^{-1}(b)$. Therefore, for any
$b_{1}, b_{2}\in B$, we have
\begin{align*}
&\quad(\chi-\chi^{\beta, \gamma})(b_{1}, b_{2})\\
&=\chi(b_{1}, b_{2})-\beta\circ\chi\Big(\gamma^{-1}(b_{1}), \gamma^{-1}(b_{2})\Big)\\
&=\varsigma(b_{1})\varsigma(b_{2})-\varsigma(b_{1}b_{2})
-\beta\Big(\varsigma(\gamma^{-1}(b_{1}))\varsigma(\gamma^{-1}(b_{2}))
-\varsigma(\gamma^{-1}(b_{1})\gamma^{-1}(b_{2}))\Big)\\
&=\varsigma(b_{1})(\varsigma-\alpha\circ\varsigma\circ\gamma^{-1})(b_{2})
+(\varsigma-\alpha\circ\varsigma\circ\gamma^{-1})(b_{1})\varsigma(b_{2})\\
&\quad-(\varsigma-\alpha\circ\varsigma\circ\gamma^{-1})(b_{1})
(\varsigma-\alpha\circ\varsigma\circ\gamma^{-1})(b_{2})
-(\varsigma-\beta\circ\varsigma\circ\gamma^{-1})(b_{1}b_{2})\\
&=b_{1}\tr\vartheta(b_{2})+\vartheta(b_{1})\tl b_{2}
-\vartheta(b_{1})\vartheta(b_{2})-\eta(b_{1}b_{2}).
\end{align*}
Similarly, one can check that $(\tr-\tr^{\beta, \gamma})(b, a)=\vartheta(b)a$ and
$(\tl-\tl^{\beta, \gamma})(a, b)=a\vartheta(b)$, for any $a\in A$ and $b\in B$.
Thus, it follows from Definition \ref{def:coho} that the non-abelian $2$-cocycles
$(\tr, \tl, \chi)$ and $(\tr^{\beta, \gamma}, \tl^{\beta, \gamma},
\chi^{\beta, \gamma})$ are equivalent.

Conversely, if the non-abelian $2$-cocycles $(\tr, \tl, \chi)$ and
$(\tr^{\beta, \gamma}, \tl^{\beta, \gamma}, \chi^{\beta, \gamma})$ are equivalent,
there exists a linear map $\lambda: B\rightarrow A$ such that Eqs.
(\ref{eq1}) and (\ref{eq2}) are satisfied, and the diagram:
$$
\xymatrixrowsep{0.5cm}
\xymatrix{0\ar[r]& A\ar[r]^{\iota}\ar@{=}[d]& E\ar[r]^{\pi}\ar[d]^{\varphi}
& B\ar[r]\ar@{=}[d]&0\\
0\ar[r]& A\ar[r]^{\iota_{A}\qquad\;} & A\rtimes_{(\tr^{\beta, \gamma},
\tl^{\beta, \gamma}, \chi^{\beta, \gamma})}B \ar[r]^{\qquad\;\pi_{B}}& B\ar[r]&0 }
$$
is a commutative in the category of perm algebras. Since $A\oplus B\rightarrow E$,
$(a, b)\mapsto a+\varsigma(b)$ is an isomorphism of vector spaces, we write an element
in $E$ as the form $a+\varsigma(b)$. Then we have
$$
(a_{1}+\varsigma(b_{1}))\diamond(a_{2}+\varsigma(b_{2}))=
a_{1}a_{2}+b_{1}\tr a_{2}+a_{1}\tl b_{2}+\chi(b_{1}, b_{2})+\varsigma(b_{1}b_{2}),
$$
for any $a_{1}, a_{2}\in A$, $b_{1}, b_{2}\in B$. For the pair of automorphisms
$(\beta, \gamma)\in\Aut(A)\times\Aut(B)$, we define a linear map $\alpha: E\rightarrow E$ by
$$
\alpha(a+\varsigma(b))=\beta(a)-\lambda(\gamma(b))+\varsigma(\gamma(b)),
$$
for any $a\in A$ and $b\in B$. Since $\beta$ and $\gamma$ are isomorphisms,
it is easy to see that $\alpha$ is an isomorphism of vector spaces. Moreover, note that
\begin{align*}
&\; \alpha(a_{1}+\varsigma(b_{1}))\diamond\alpha(a_{2}+\varsigma(b_{2}))\\
=&\; \Big(\beta(a_{1})-\lambda(\gamma(b_{1}))+\varsigma(\gamma(b_{1}))\Big)
\Big(\beta(a_{2})-\lambda(\gamma(b_{2}))+\varsigma(\gamma(b_{2}))\Big)\\
=&\; \beta(a_{1})\beta(a_{2})+\gamma(b_{1})\tr^{\beta,\gamma}\beta(a_{2})
+\beta(a_{1})\tl^{\beta,\gamma}\gamma(b_{2})+\chi^{\beta,\gamma}(\gamma(b_{1}),
\gamma(b_{2}))\\[-1mm]
&\qquad\qquad\qquad\qquad\qquad\qquad\qquad+\lambda(\gamma(b_{1})\gamma(b_{2}))
+\varsigma(\gamma(b_{1})\gamma(b_{2}))\\
=&\;\beta\Big(a_{1}a_{2}+b_{1}\tr a_{2}+a_{1}\tl b_{2}+\chi(b_{1}, b_{2})\Big)
-\lambda(\gamma(b_{1}b_{2}))+\varsigma(\gamma(b_{1}b_{2}))\\
=&\; \alpha\Big(a_{1}a_{2}+b_{1}\tr a_{2}+a_{1}\tl b_{2}
+\chi(b_{1}, b_{2})+\varsigma(b_{1}b_{2})\Big)\\
=&\; \alpha\Big((a_{1}+\varsigma(b_{1}))\diamond(a_{2}+\varsigma(b_{2}))\Big),
\end{align*}
for any $a_{1}, a_{2}\in A$ and $b_{1}, b_{2}\in B$, we get $\alpha\in\Aut(E)$.
Clearly, $\alpha|_{A}=\beta$ and $\bar{\alpha}=\pi\circ\alpha\circ\varsigma=\gamma$.
Hence the pair $(\beta, \gamma)$ is inducible.
\end{proof}

From Eqs. (\ref{eq1}) and (\ref{eq2}), we can directly get the following corollary.

\begin{cor}\label{cor:describe}
Let $\mathcal{E}: \xymatrix@C=0.5cm{0\ar[r]& A\ar[r]^{\iota}& E\ar[r]^{\pi}& B\ar[r]&0}$
be a non-abelian extension of perm algebra $B$ by perm algebra $A$, and $(\tr, \tl, \chi)$
be the corresponding non-abelian $2$-cocycle of $\mathcal{E}$ with section $\varsigma$.
A pair $(\beta, \gamma)\in\Aut(A)\times\Aut(B)$ of automorphisms is inducible if and only
if there exists a linear map $\lambda: B\rightarrow A$ such that for any
$a\in A$ and $b, b_{1}, b_{2}\in B$,
\begin{align*}
\beta(b\tr a)&=\gamma(b)\tr\beta(a)-\lambda(\gamma(b))\beta(a),\\
\beta(a\tl b)&=\beta(a)\tl\gamma(b)-\beta(a)\lambda(\gamma(b)),\\
\beta(\chi(b_{1}, b_{2}))&=\chi(\gamma(b_{1}), \gamma(b_{2}))
-\gamma(b_{1})\tr\lambda(\gamma(b_{2}))+\lambda(\gamma(b_{1})\gamma(b_{2}))\\
&\quad-\lambda(\gamma(b_{1}))\tl\gamma(b_{2})+\lambda(\gamma(b_{1}))
\lambda(\gamma(b_{2})).
\end{align*}
\end{cor}

For a non-abelian extension $\mathcal{E}: \xymatrix@C=0.5cm{0\ar[r]& A\ar[r]^{\iota}
& E\ar[r]^{\pi}& B\ar[r]&0}$ of perm algebras, we define
$$
\Aut_{A}^{\id}(E):=\{\alpha\in\Aut_{A}(E) \mid\kappa(\alpha)=(\id_{A}, \id_{B})\}.
$$
Then $\Aut_{A}^{\id}(E)$ is a subgroup of $\Aut_{A}(E)$, and by the definition
of map $\kappa$, there exists an exact sequence as follows.

\begin{cor}\label{cor: exa}
Let $\mathcal{E}: \xymatrix@C=0.5cm{0\ar[r]& A\ar[r]^{\iota} & E\ar[r]^{\pi}& B\ar[r]&0}$
be a non-abelian extension of perm algebras. With the above notations, there is an
exact sequence of groups
$$
\xymatrix{1\ar[r]& \Aut_{A}^{\id}(E)\ar[r]^{\sigma} & \Aut_{A}(E)\ar[r]^{\kappa\qquad}
& \Aut(A)\times\Aut(B)},
$$
where $\sigma$ is the inclusion map.
\end{cor}

\subsection{Fundamental sequence of Wells} \label{subsec:well}
We interpret the above theorem in terms of the Wells map in the context of
perm algebras. Let $\mathcal{E}: \xymatrix@C=0.5cm{0\ar[r]& A\ar[r]^{\iota}
& E\ar[r]^{\pi}& B\ar[r]&0}$ be a non-abelian extension of perm algebra $B$ by $A$,
with section $\varsigma$. If the multiplication of $A$ is trivial, i.e., $a_{1}a_{2}=0$ for
any $a_{1}, a_{2}\in A$, then $\mathcal{E}$ is an abelian extension, and there is a 
$B$-bimodule structure on $A$ by the multiplication $``\diamond"$ of $E$, i.e.,
$b\tr a=\varsigma(b)\diamond a$ and $a\tr b=a\diamond\varsigma(b)$ for $a\in A$, $b\in B$.
Moreover, one can check that this bimodule does not depend on the choice of $\varsigma$.
We fix the $B$-bimodule structure on $A$ as above. Then the corresponding non-abelian
$2$-cocycle $(\tr, \tl, \chi)$ of $\mathcal{E}$ is given by the bimodule structure and
a $2$-cocycle $\chi$ of perm algebra $B$ with coefficients in bimodule $A$, and
$H^{2}_{nab}(B, A)$ is isomorphic to the second cohomology group of $B$ with
coefficients in bimodule $A$. Thus, in this case, $H^{2}_{nab}(B, A)$ is an additive group,
and for convenience, we record the elements in $H^{2}_{nab}(B, A)$ by $[\chi]$.

Now let
$$
\mathcal{E}:\quad \xymatrix{0\ar[r]& A\ar[r]^{\iota}& E\ar[r]^{\pi}& B\ar[r]&0}
$$
be an abelian extension of perm algebra $B$ by $A$, where the $B$-bimodule structure on
$A$ is given by $(\tr, \tl)$. Denote by $\Aut(A)$ the set of all linear automorphisms on $A$,
and by $\Aut_{\tr, \tl}(A, B)$ the set of all pairs $(\beta, \gamma)\in\Aut(A)
\times\Aut(B)$ satisfy
$$
\beta(b\tr a)=\gamma(b)\tr\beta(a),\qquad\qquad   \beta(a\tl b)=\beta(a)\tl\gamma(b),
$$
for $a\in A$ and $b\in B$. Then $\Aut_{\tr, \tl}(A, B)$ is a subgroup of
$\Aut(A)\times\Aut(B)$. Let $\varsigma$ be a section of $\mathcal{E}$, and $\chi$ be
the 2-cocycle corresponding to $\mathcal{E}$ with $\varsigma$. For a pair
$(\beta, \gamma)\in\Aut(B)\times\Aut(A)$, we define $\chi^{\beta, \gamma}: B\times B
\rightarrow A$ by
$$
\chi^{\beta, \gamma}(b_{1}, b_{2})=\beta\circ\chi(\gamma^{-1}(b_{1}), \gamma^{-1}(b_{2})),
$$
for $b_{1}, b_{2}\in B$. If $(\beta, \gamma)\in\Aut_{\tr, \tl}(A, B)$,
by direct calculation, one can check that $(\tr, \tl, \chi^{\beta, \gamma})\in
Z^{2}_{nab}(B, A)$. Following from Theorem \ref{thm: induc} and Corollary
\ref{cor:describe}, we obtain the following Proposition.

\begin{pro}\label{pro: abe-induc}
Let $\mathcal{E}: \xymatrix@C=0.5cm{0\ar[r]& A\ar[r]^{\iota}& E\ar[r]^{\pi}& B\ar[r]&0}$
be an abelian extension of perm algebra $B$ by bimodule $(A, \tr, \tl)$,
$\varsigma$ be a section of $\mathcal{E}$, and $\chi$ be the 2-cocycle corresponding
to $\mathcal{E}$ with $\varsigma$. Then a pair $(\beta, \gamma)\in\Aut_{\tr, \tl}(A, B)$
is inducible if and only if the 2-cocycles $\chi$ and $\chi^{\beta, \gamma}$ are equivalent,
if and only if there exists a linear map $\lambda: B\rightarrow A$ such that for any
$b_{1}, b_{2}\in B$,
\begin{align*}
\beta(\chi(b_{1}, b_{2}))&=\chi(\gamma(b_{1}),
\gamma(b_{2}))-\gamma(b_{1})\tr\lambda(\gamma(b_{2}))\\
&\qquad+\lambda(\gamma(b_{1})\gamma(b_{2}))-\lambda(\gamma(b_{1}))\tl\gamma(b_{2}).
\end{align*}
\end{pro}

For an abelian extension $\mathcal{E}: \xymatrix@C=0.5cm{0\ar[r]& A\ar[r]^{\iota}
& E\ar[r]^{\pi}& B\ar[r]&0}$, we define the Wells map by
$$
\mathcal{W}_{\mathcal{E}}:\; \Aut_{\tr, \tl}(A, B)\rightarrow \HH_{nab}^{2}(B, A)
\qquad\qquad (\beta, \gamma)\mapsto[\chi^{\beta, \gamma}-\chi].
$$
where $\chi$ be the 2-cocycle corresponding to $\mathcal{E}$. The map
$\mathcal{W}_{\mathcal{E}}$ is well-defined and a group homomorphism since
$\HH_{nab}^{2}(B, A)$ is an additive group in this case. Now, we get the 
fundamental sequence of Wells in the following form.

\begin{pro}[Fundamental sequence of Wells]\label{pro:wells}
Let $\mathcal{E}: \xymatrix@C=0.5cm{0\ar[r]& A\ar[r]^{\iota} & E\ar[r]^{\pi}& B\ar[r]&0}$
be an abelian extension of perm algebra $B$ by bimodule $A$.
With the above notations, there is an exact sequence of groups:
$$
\xymatrix{1\ar[r]&\Aut_{A}^{\id}(E)\ar[r]^{\sigma}&\Aut_{A}(E)\ar[r]^{\kappa \quad }
&\Aut_{\tr, \tl}(A, B)\ar[r]^{\; \mathcal{W}_{\mathcal{E}}} & \HH_{nab}^{2}(B, A)},
$$
where $\sigma$ is the inclusion map.
\end{pro}

\begin{proof}
First, we need to show it is exact at $\Aut_{A}(E)$. Note that an element $\alpha\in
\Aut_{A}(E)$ such that $\alpha\in\Ker(\kappa)$ if and only if $\kappa(\alpha)=
(\id_{A}, \id_{B})$, i.e., $\alpha\in\Aut_{A}^{\id}(E)$, we obtain the sequence is exact 
at $\Aut_{A}(E)$. Second, by Proposition \ref{pro: abe-induc}, we get a pair $(\beta, \gamma)
\in\Aut_{\tr, \tl}(A, B)$ satisfying $(\beta, \gamma)\in\Img(\kappa)$, i.e.,
$(\beta, \gamma)$ is inducible, if and only if $\mathcal{W}_{\mathcal{E}}(\beta, \gamma)=0$.
This means that the sequence is exact at $\Aut_{\tr, \tl}(A, B)$. And so that,
the sequence is exact in the category of groups.
\end{proof}

In particular, if the abelian extension $\mathcal{E}: \xymatrix@C=0.5cm{0\ar[r]
& A\ar[r]^{\iota} & E\ar[r]^{\pi}& B\ar[r]&0}$ is split in the category
of perm algebras, then the perm algebra $E$ is isomorphic to the
perm algebra $A\rtimes B$. Let $\varsigma: B\rightarrow E$, $b\mapsto(a, b)$. Then
$\varsigma$ is a section of $\mathcal{E}$. In this case, the corresponding 2-cocycle
$\chi$ of $\mathcal{E}$ with section $\varsigma$ is zero in $\HH_{nab}^{2}(B, A)$,
and the Wells map vanishes identically. Thus we have an exact sequence in the category
of groups:
$$
\xymatrix{1\ar[r]&\Aut_{A}^{\id}(E)\ar[r]^{\sigma}& \Aut_{A}(E)\ar[r]^{\kappa\quad}
& \Aut_{\tr, \tl}(A, B)\ar[r]&1}.
$$
We define map $\varrho: \Aut_{\tr, \tl}(A, B)\rightarrow\Aut_{A}(E)$ by
$$
\varrho(\beta, \gamma)(a, b)=(\beta(a), \gamma(b)),
$$
for any $(\beta, \gamma)\in\Aut_{\tr, \tl}(A, B)$ and $(a, b)\in A\oplus B$.
Then $\kappa(\varrho(\beta, \gamma))=(\varrho(\beta, \gamma)|_{A},\;
\overline{\varrho(\beta, \gamma)})=(\beta, \gamma)$, since $\pi\circ\varrho
(\beta, \gamma)\circ\varsigma(b)=\gamma(b)$ for any $b\in B$.
That is to say, the exact sequence above is split in the category of groups.
Thus we have the following proposition.

\begin{pro}\label{pro:split}
Let $\mathcal{E}: \xymatrix@C=0.5cm{0\ar[r] & A\ar[r]^{\iota} & E\ar[r]^{\pi}& B\ar[r]&0}$
be a split abelian extension of perm algebra $B$ by bimodule $A$. Then
$$
\Aut_{A}(E)\cong \Aut_{A}^{\id}(E)\rtimes\Aut_{\tr, \tl}(A, B)
$$
as groups, where $\Aut_{A}^{\id}(E)\rtimes\Aut_{\tr, \tl}(A, B)$ is the semi-direct product
of groups.
\end{pro}

To finish this section, we present a simple example.

\begin{ex}\label{ex:auto}
Let $B=\mathbb{C}\{e_{1}, e_{2}\}$ be a 2-dimensional perm algebra with non-zero products:
$e_{1}e_{1}=e_{1}$ and $e_{2}e_{1}=e_{2}$, and $A=\mathbb{C}\{e_{3}\}$ be a
$B$-bimodule under the module action: $e_{1}\tr e_{3}=e_{3}=e_{3}\tl e_{1}$
and $e_{2}\tr e_{3}=0=e_{3}\tl e_{2}$. Then we get a split abelian extension
$\xymatrix@C=0.5cm{0\ar[r] & A\ar[r]^{\iota} & E\ar[r]^{\pi}& B\ar[r]&0}$,
where $E=\mathbb{C}\{e_{1}, e_{2}, e_{3}\}$ with non-zero products: $e_{1}e_{1}=e_{1}$,
$e_{2}e_{1}=e_{2}$ and $e_{1}e_{3}=e_{3}=e_{3}e_{1}$.
By direct computations, we can obtain $\Aut(A)=\bk^{\ast}:=\bk\backslash\{0\}$,
$$
\Aut(B)=\left\{{\tiny\left(\begin{array}{ccc} 1&0\\ a&b
\end{array}\right)}\;\Big|\; a\in\bk, b\in\bk^{\ast}\right\},\qquad
\Aut(E)=\left\{{\tiny\left(\begin{array}{ccc} 1&0&0\\ a&b&0\\ 0&0&c
\end{array}\right)}\;\Big|\; a\in\bk, b, c\in\bk^{\ast}\right\},
$$
$\Aut_{\tr, \tl}(A, B)=\Aut(A)\times\Aut(B)$ and $\Aut_{A}^{\id}(E)=1$.
Then, it is easy to see that $\Aut_{A}(E)\cong\Aut(E)\cong\Aut(A)\times\Aut(B)\cong
\Aut_{A}^{\id}(E)\rtimes\Aut_{\tr, \tl}(A, B)$.
\end{ex}

\section{Perm bialgebras and $\mathcal{S}$-equation in perm algebas} \label{sec: YB-equ}

In this section, we consider a class of special matched pairs of perm algebras.
We introduce the notions of Manin triple for perm algebras and perm bialgebras.
Under certain conditions, we show that perm bialgebras, standard Manin triples
for perm algebras and certain matched pairs of perm algebras are equivalent.
We introduce and study coboundary perm bialgebras and our study leads to the
''$\mathcal{S}$-equation" in perm algebras, which is an analogue of the classical
Yang-Baxter equation. A symmetric solution of $\mathcal{S}$-equation gives a perm
bialgebra.

From now on, all vector spaces, algebras and modules are finite-dimensional over
field $\bk$. For any vector space $V$, we denote $V^{\ast}$ the dual space of $V$.
Let $V$, $W$ be two $\bk$ vector spaces. for any $v\in V$ and $u\in V^{\ast}$,
we denote $\langle u, v\rangle:=u(v)\in\bk$. For a linear map $f: V\rightarrow W$,
we define the map $f^{\ast}: W^{\ast}\rightarrow V^{\ast}$ by $\langle f^{\ast}(w), v\rangle
=\langle w, f(v)\rangle$ for any $v\in V$ and $w\in W^{\ast}$.

Let $A$ and $B$ be two perm algebras, $\kl_{1}, \kr_{1}: A\rightarrow\End_{\bk}(B)$
and $\kl_{2}, \kr_{2}: B\rightarrow\End_{\bk}(A)$ be four linear maps.
Define bilinear linear maps
\begin{align*}
\br: A\times B\rightarrow B &\qquad\quad a\br b=\kl_{1}(a)(b),\\
\bl: B\times A\rightarrow B &\qquad\quad b\bl a=\kr_{1}(a)(b),\\
\tr: B\times A\rightarrow A &\qquad\quad b\tr a=\kl_{2}(b)(a),\\
\tl: A\times B\rightarrow A &\qquad\quad a\tl b=\kr_{2}(b)(a),
\end{align*}
for any $a\in A$ and $b\in B$. It is easy to see that $(B, \br, \bl)$ is a bimodule
over $A$ if and only if $(B, \kl_{1}, \kr_{1})$ is a representation of $A$, and
$(A, \tr, \tl)$ is a bimodule over $B$ if and only if $(A, \kl_{2}, \kr_{2})$ is a
representation of $B$. We also call $(A, B, \kl_{1}, \kr_{1}, \kl_{2}, \kr_{2})$
is a matched pair if $(A, B, \br, \bl, \tr, \tl)$ is a matched pair.
Then, by the definition of matched pair $(A, B, \br, \bl, \tr, \tl)$, we have the
following proposition.

\begin{pro}\label{pro:rep-match}
Let $A$, $B$ be two perm algebras, $(B, \kl_{1}, \kr_{1})$ be a representation of $A$,
and $(A, \kl_{2}, \kr_{2})$ be a representation of $B$ respectively. Then $(A, B, \kl_{1},
\kr_{1}, \kl_{2}, \kr_{2})$ is a matched pair if and only if,
for any $a, a_{1}, a_{2}\in A$, $b, b_{1}, b_{2}\in B$,
\begin{align}
\kl_{2}(b)(a_{1}a_{2})&=\kl_{2}(b)(a_{2}a_{1})=(\kl_{2}(b)(a_{1}))a_{2}
+\kl_{2}(\kr_{1}(a_{1})(b))(a_{2}),                                            \label{mat1}\\
\kr_{2}(b)(a_{1}a_{2})&=a_{1}(\kr_{2}(b)(a_{2}))
+\kr_{2}(\kl_{1}(a_{2})(b))(a_{1})                                             \label{mat2}\\
&=a_{1}(\kl_{2}(b)(a_{2}))+\kr_{1}(\kr_{1}(a_{2})(b))(a_{1})
=(\kr_{2}(b)(a_{1}))a_{2}+\kl_{2}(\kl_{1}(a_{1})(b))(a_{2}),                   \nonumber\\
\kl_{1}(a)(b_{1}b_{2})&=\kl_{1}(a)(b_{2}b_{1})
=(\kl_{1}(a)(b_{1}))b_{2}+\kl_{1}(\kr_{2}(b_{1})(a))(b_{2}),                   \label{mat3}\\
\kr_{1}(a)(b_{1}b_{2})&=b_{1}(\kr_{1}(a)(b_{2}))
+\kr_{1}(\kl_{2}(b_{2})(a))(b_{1})                                             \label{mat4}\\
&=b_{1}(\kl_{1}(a)(b_{2}))+\kr_{1}(\kr_{2}(b_{2})(a))(b_{1})
=(\kr_{1}(a)(b_{1}))b_{2}+\kl_{1}(\kl_{2}(b_{1})(a))(b_{2}).                   \nonumber
\end{align}
\end{pro}

Define a bilinear operation $\ast: (A\oplus B)\times(A\oplus B)\rightarrow A\oplus B$ by
$$
(a_{1}, b_{1})\ast(a_{2}, b_{2})=\Big(a_{1}a_{2}+\kl_{2}(b_{1})(a_{2})
+\kr_{2}(b_{2})(a_{1}),\ \
b_{1}b_{2}+\kl_{1}(a_{1})(b_{2})+\kr_{1}(a_{2})(b_{1})\Big),
$$
where $a_{1}, a_{2}\in A$ and $b_{1}, b_{2}\in B$. Then by Proposition \ref{Pro:exten},
$(A\oplus B, \ast)$ is a prem algebra if and only if $(A, B, \kl_{1}, \kr_{1}, \kl_{2},
\kr_{2})$ is a matched pair of prem algebras. The prem algebra $(A\oplus B, \ast)$ is
denoted by $A\bowtie^{\kl_{1}, \kr_{1}}_{\kl_{2}, \kr_{2}}B$, or simply $A\bowtie B$,
and $(A, B, \kl_{1}, \kr_{1}, \kl_{2}, \kr_{2})$ is said to be the corresponding matched
pair of $A\bowtie B$.

\begin{lem}\label{lem:dual}
Let $A$ be a perm algebra and $(V, \kl, \kr)$ be a representation of $A$. Define linear maps
\begin{align*}
& \kl^{\ast}: A\rightarrow\End_{\bk}(V^{\ast}), \qquad\quad
\langle\kl^{\ast}(a)(u),\, v\rangle=\langle u,\, \kl(a)(v)\rangle,\\
& \kr^{\ast}: A\rightarrow\End_{\bk}(V^{\ast}), \qquad\quad
\langle\kr^{\ast}(a)(u),\, v\rangle=\langle u,\, \kr(a)(v)\rangle,
\end{align*}
for any $a\in A$, $v\in V$ and $u\in V^{\ast}$. Then $(V^{\ast}, \kr^{\ast}-\kl^{\ast},
\kr^{\ast})$ is also a representation of $A$. We call it the dual representation of
$(V, \kl, \kr)$.
\end{lem}

\begin{proof}
For any $a_{1}, a_{2}\in A$, $v\in V$ and $u\in V^{\ast}$, we have
\begin{align*}
\langle(\kr^{\ast}-\kl^{\ast})(a_{1}a_{2})(u),\; v\rangle
&=\langle u,\; \kr(a_{1}a_{2})(v)\rangle-\langle u,\; \kl(a_{1}a_{2})(v)\rangle\\
\langle((\kr^{\ast}-\kl^{\ast})(a_{1})\circ(\kr^{\ast}-\kl^{\ast})(a_{2}))(u),\; v\rangle
&=\langle(\kr^{\ast}-\kl^{\ast})(a_{2})(u),\; \kr(a_{1})(v)-\kl(a_{1})(v)\rangle\\
&=\langle u,\; (\kr(a_{2})\circ\kr(a_{1}))(v)\rangle
-\langle u,\; (\kl(a_{2})\circ\kr(a_{1}))(v)\rangle\\
&\qquad-\langle u,\; (\kr(a_{2})\circ\kl(a_{1}))(v)\rangle
+\langle u,\; (\kl(a_{2})\circ\kl(a_{1}))(v)\rangle,\\
&=\langle u,\; (\kr(a_{2})\circ\kr(a_{1}))(v)\rangle
-\langle u,\; (\kr(a_{2})\circ\kl(a_{1}))(v)\rangle
\end{align*}
\begin{align*}
\langle((\kr^{\ast}-\kl^{\ast})(a_{1})\circ\kr^{\ast}(a_{2}))(u),\; v\rangle
&=\langle(\kr^{\ast}(a_{2})(u),\; \kr(a_{1})(v)-\kl(a_{1})(v)\rangle\\
&=\langle u,\; (\kr(a_{2})\circ\kr(a_{1}))(v)\rangle
-\langle u,\; (\kr(a_{2})\circ\kl(a_{1}))(v)\rangle,\\
\langle(\kr^{\ast}(a_{2})\circ(\kr^{\ast}-\kl^{\ast})(a_{1}))(u),\; v\rangle
&=\langle((\kr^{\ast}-\kl^{\ast})(a_{1}))(u),\; \kr(a_{2})(v)\rangle\\
&=\langle u,\; (\kr(a_{1})\circ\kr(a_{2}))(v)\rangle
-\langle u,\; (\kl(a_{1})\circ\kr(a_{2}))(v)\rangle.
\end{align*}
Since $(V, \kl, \kr)$ be a representation of $A$, we get $(\kr^{\ast}-\kl^{\ast})(a_{1}a_{2})
=(\kr^{\ast}-\kl^{\ast})(a_{1})\circ(\kr^{\ast}-\kl^{\ast})(a_{2})=(\kr^{\ast}
-\kl^{\ast})(a_{1})\circ\kr^{\ast}(a_{2})=\kr^{\ast}(a_{2})\circ(\kr^{\ast}-\kl^{\ast})(a_{1})$.
Moreover, it is easy to see that
$\kr^{\ast}(a_{1}a_{2})=\kr^{\ast}(a_{2})\circ\kr^{\ast}(a_{1})
=\kr^{\ast}(a_{1})\circ\kr^{\ast}(a_{2})$, for any $a_{1}, a_{2}\in A$. Thus $(V^{\ast},
\kr^{\ast}-\kl^{\ast}, \kr^{\ast})$ is a representation of $A$.
\end{proof}

\subsection{Manin triple for perm algebras and perm bialgebras} \label{subsec:bialg}

For a perm algebra $A$, we consider the regular representation $(A, \kl_{A}, \kr_{A})$.
Then we get a dual representation $(A^{\ast}, \kr_{A}^{\ast}-\kl_{A}^{\ast}, \kr_{A}^{\ast})$.
For the relation of this two representations, we need consider some special bilinear
form on $A$. A bilinear form $\omega$ on a perm algebra $A$ is called {\it invariant}
if for any $a_{1}, a_{2}, a_{3}\in A$,
$$
\omega(a_{1}a_{2}, a_{3})=\omega(a_{2}, a_{3}a_{1})-\omega(a_{2}, a_{1}a_{3}).
$$

\begin{pro}\label{pro:rep-iso}
Let $A$ be a perm algebra. If there exists a nondegenerate skew-symmetric invariant
bilinear form on $A$, two representations $(A, \kl_{A}, \kr_{A})$ and
$(A^{\ast}, \kr_{A}^{\ast}-\kl_{A}^{\ast}, \kr_{A}^{\ast})$ are isomorphic.
If the two representations $(A, \kl_{A}, \kr_{A})$ and
$(A^{\ast}, \kr_{A}^{\ast}-\kl_{A}^{\ast}, \kr_{A}^{\ast})$ are isomorphic,
there exists a nondegenerate invariant bilinear form on $A$.
\end{pro}

\begin{proof}
If $\omega$ is a nondegenerate skew-symmetric invariant bilinear form on $A$, we define
linear map $\varphi: A\rightarrow A^{\ast}$ by $\langle \varphi(a_{1}), a_{2}\rangle=
\omega(a_{1}, a_{2})$ for all $a_{1}, a_{2}\in A$. Since $\omega$ is nondegenerate, $\varphi$
is a bijection. Note that
\begin{align*}
\langle\varphi(\kl_{A}(a_{1})a_{2}),\; a_{3}\rangle&=\omega(a_{1}a_{2},\; a_{3})
=\omega(a_{2}, a_{3}a_{1})-\omega(a_{2}, a_{1}a_{3})\\
&=\langle\kr^{\ast}_{A}(a_{1})\varphi(a_{2}),\; a_{3}\rangle
-\langle\kl^{\ast}_{A}(a_{1})\varphi(a_{2}),\; a_{3}\rangle,\\
\langle\varphi(\kr_{A}(a_{1})a_{2}),\; a_{3}\rangle&=\omega(a_{2}a_{1},\; a_{3})
=\omega(a_{1}, a_{3}a_{2})-\omega(a_{1}, a_{2}a_{3})\\
&=\omega(a_{2},\; a_{3}a_{1})-\omega(a_{2},\; a_{1}a_{3})
+\omega(a_{3},\; a_{1}a_{2})-\omega(a_{3},\; a_{2}a_{1})\\
&=\langle\kr^{\ast}_{A}(a_{1})\varphi(a_{2}),\; a_{3}\rangle,
\end{align*}
for any $a_{1}, a_{2}, a_{3}\in A$, we obtain that $\varphi$ is an isomorphism of
representation of $A$. The converse can be proved similarly. We omit the details.
\end{proof}

We introduce the notion of Manin triple for perm algebras.

\begin{defi} \label{Def:Manin-tri}
A {\rm Manin triple} of perm algebras is a triple of perm algebras $(\mathcal{A},
\mathcal{A}^{+}, \mathcal{A}^{-})$ together with a nondegenerate skew-symmetric invariant
bilinear form $\omega$ on $\mathcal{A}$ such that the following conditions are satisfied:
\begin{enumerate}
\item[$(i)$] $\mathcal{A}^{+}$, $\mathcal{A}^{-}$ are perm subalgebras of $\mathcal{A}$;

\item[$(ii)$] $\mathcal{A}=\mathcal{A}^{+}\oplus\mathcal{A}^{-}$ as vector spaces;

\item[$(iii)$] $\mathcal{A}^{+}$ and $\mathcal{A}^{-}$ are isotropic with respect
     to $\omega$, that is, $\omega(a_{1}^{+}, a_{2}^{+})=0=\omega(a_{1}^{-}, a_{2}^{-})$
     for all $a_{1}^{+}, a_{2}^{+}\in\mathcal{A}^{+}$ and
     $a_{1}^{-}, a_{2}^{-}\in\mathcal{A}^{-}$.
\end{enumerate}
\end{defi}

A {\it morphism between two Manin triples} of perm algebras $(\mathcal{A}, \mathcal{A}^{+},
\mathcal{A}^{-})$ and $(\mathcal{B}, \mathcal{B}^{+}, \mathcal{B}^{-})$
associated to two nondegenerate skew-symmetric invariant bilinear forms $\omega_{1}$ and
$\omega_{2}$ respectively, is a morphism of prem algebras $f: \mathcal{A}\rightarrow
\mathcal{B}$ such that for any $a_{1}, a_{2}\in A$,
$$
f(\mathcal{A}^{+})\subseteq\mathcal{B}^{+},\qquad
f(\mathcal{A}^{-})\subseteq\mathcal{B}^{-},\qquad
\omega_{1}(a_{1}, a_{2})=\omega_{2}(f(a_{1}), f(a_{2})).
$$
If in addition, $f$ is a bijection, then the two Manin triples are called {\it isomorphic}.
Let $A$ be a perm algebra. If
\begin{enumerate}
\item[(a)] there exist a perm algebra structure on the dual space $A^{\ast}$;

\item[(b)] there exist a perm algebra structure on the direct sum $A\oplus A^{\ast}$
     such that $A$ and $A^{\ast}$ are perm subalgebras of $A\oplus A^{\ast}$;

\item[(c)] and the bilinear form on $A\oplus A^{\ast}$:
     $$
     \widehat{\omega}((a_{1}, b_{1}),\, (a_{2}, b_{2}))=\langle a_{1}, b_{2}\rangle
     -\langle a_{2}, b_{1}\rangle
     $$
     for any $a_{1}, a_{2}\in A$ and $b_{1}, b_{2}\in A^{\ast}$, is invariant,
\end{enumerate}
then $(A\oplus A^{\ast}, A, A^{\ast})$ is called a {\it (standard) Manin triple of perm
algebras associated to $\widehat{\omega}$}. Obviously, a standard Manin triple of perm
algebras is a Manin triple of perm algebras. Conversely, we have

\begin{pro}\label{pro:main-iso}
Every Manin triple of perm algebras is isomorphic to a standard one.
\end{pro}

\begin{proof}
Let perm algebras $(\mathcal{A}, \mathcal{A}^{+}, \mathcal{A}^{-})$ together with a
nondegenerate skew-symmetric invariant bilinear form $\omega$ be a Manin triple.
Since $\mathcal{A}^{+}$ and $\mathcal{A}^{-}$ are isotropic under the nondegenrate
bilinear form $\omega$, $(\mathcal{A}^{+})^{\ast}$ and $\mathcal{A}^{-}$ are identified
by $\omega$, and the perm algebra structure on $\mathcal{A}^{-}$ is transferred to
$(\mathcal{A}^{+})^{\ast}$. And so that, there exist a perm algebra structure on
$\mathcal{A}^{+}\oplus(\mathcal{A}^{+})^{\ast}$. Hence, we can transfer the bilinear form
$\omega$ to $\mathcal{A}^{+}\oplus(\mathcal{A}^{+})^{\ast}$, we obtain the standard
bilinear form $\widehat{\omega}$ defined as above. Thus, the Manin triple $(\mathcal{A},
\mathcal{A}^{+}, \mathcal{A}^{-})$ is isomorphic to the standard Manin triple
$(\mathcal{A}^{+}\oplus(\mathcal{A}^{+})^{\ast}, \mathcal{A}^{+}, (\mathcal{A}^{+})^{\ast})$.
\end{proof}

In fact, there is a close relationship between Manin triple and matched pair of perm algebras.
Let $A$ be a perm algebra. If there is a perm algebra structure on $A^{\ast}$, we denote
$(A^{\ast}, \kl_{A^{\ast}}, \kr_{A^{\ast}})$ the regular representation of $A^{\ast}$. Then
$(A, \kr_{A^{\ast}}^{\ast}-\kl_{A^{\ast}}^{\ast}, \kr_{A^{\ast}}^{\ast})$ is also a
representation of $A^{\ast}$ by Lemma \ref{lem:dual}.

\begin{pro}\label{pro:main-mat}
Let $A$ be a perm algebra. Suppose that there is a perm algebra structure on $A^{\ast}$.
Then there exists a perm algebra structure on the vector space $A\oplus A^{\ast}$
such that $(A\oplus A^{\ast}, A, A^{\ast})$ is a standard Manin triple with respect
to $\widehat{\omega}$ defined as above if and only if $(A, A^{\ast},
\kr_{A}^{\ast}-\kl_{A}^{\ast}, \kr_{A}^{\ast}, \kr_{A^{\ast}}^{\ast}-\kl_{A^{\ast}}^{\ast},
\kr_{A^{\ast}}^{\ast})$ is a matched pair of perm algebras.
\end{pro}

\begin{proof}
If $(A, A^{\ast}, \kr_{A}^{\ast}-\kl_{A}^{\ast}, \kr_{A}^{\ast},
\kr_{A^{\ast}}^{\ast}-\kl_{A^{\ast}}^{\ast}, \kr_{A^{\ast}}^{\ast})$ is a matched pair of
perm algebras, then it is straightforward to show that the bilinear form $\widehat{\omega}$
is invariant on $A\bowtie A^{\ast}$. Conversely, if $(\mathcal{A}=A\oplus A^{\ast},
A, A^{\ast})$ is a standard Manin triple with respect to $\widehat{\omega}$,
then $\widehat{\omega}$ induce an isomorphism from $\mathcal{A}=A\oplus A^{\ast}$ to
$A\oplus A^{\ast}$, and the perm algebra structure on $A\oplus A^{\ast}$ just
$A\bowtie A^{\ast}$. Thus $(A, A^{\ast}, \kr_{A}^{\ast}-\kl_{A}^{\ast}, \kr_{A}^{\ast},
\kr_{A^{\ast}}^{\ast}-\kl_{A^{\ast}}^{\ast}, \kr_{A^{\ast}}^{\ast})$ is a matched pair of
perm algebras.
\end{proof}

We naturally ask when is $(A, A^{\ast}, \kr_{A}^{\ast}-\kl_{A}^{\ast}, \kr_{A}^{\ast},
\kr_{A^{\ast}}^{\ast}-\kl_{A^{\ast}}^{\ast}, \kr_{A^{\ast}}^{\ast})$ a matched pair
of perm algebras?

\begin{pro}\label{pro:mat-bia}
Let $A$ be a perm algebra. Suppose that there is a perm algebra structure on $A^{\ast}$.
Then $(A, A^{\ast}, \kr_{A}^{\ast}-\kl_{A}^{\ast}, \kr_{A}^{\ast},
\kr_{A^{\ast}}^{\ast}-\kl_{A^{\ast}}^{\ast}, \kr_{A^{\ast}}^{\ast})$ is a matched pair of
perm algebras if and only if for any $a_{1}, a_{2}\in A$, $b\in A^{\ast}$,
\begin{align}
&\kr_{A^{\ast}}^{\ast}(b)(a_{1}a_{2})-\kl_{A^{\ast}}^{\ast}(b)(a_{1}a_{2})
=\kr_{A^{\ast}}^{\ast}(b)(a_{2}a_{1})-\kl_{A^{\ast}}^{\ast}(b)(a_{2}a_{1})     \label{matc1}\\
&\qquad=\kr_{A^{\ast}}^{\ast}(b)(a_{1})a_{2}-\kl_{A^{\ast}}^{\ast}(b)(a_{1})a_{2}
+\kr_{A^{\ast}}^{\ast}(\kr_{A}^{\ast}(a_{1})(b))(a_{2})
-\kl_{A^{\ast}}^{\ast}(\kr_{A}^{\ast}(a_{1})(b))(a_{2}),                        \nonumber\\
&\kr_{A^{\ast}}^{\ast}(b)(a_{1}a_{2})=a_{1}\kr_{A^{\ast}}^{\ast}(b)(a_{2})
+\kr_{A^{\ast}}^{\ast}(\kr_{A}^{\ast}(a_{2})(b))(a_{1})
-\kr_{A^{\ast}}^{\ast}(\kl_{A}^{\ast}(a_{2})(b))(a_{1})                        \label{matc2}\\
&\qquad=a_{1}\kr_{A^{\ast}}^{\ast}(b)(a_{2})-a_{1}\kl_{A^{\ast}}^{\ast}(b)(a_{2})
+\kr_{A^{\ast}}^{\ast}(\kr_{A}^{\ast}(a_{2})(b))(a_{1})                       \nonumber\\
&\qquad=\kr_{A^{\ast}}^{\ast}(b)(a_{1})a_{2}
+\kr_{A^{\ast}}^{\ast}(\kr_{A}^{\ast}(a_{1})(b))(a_{2})
-\kr_{A^{\ast}}^{\ast}(\kl_{A}^{\ast}(a_{1})(b))(a_{2})                       \nonumber\\
&\qquad\qquad\qquad\qquad\qquad-\kl_{A^{\ast}}^{\ast}(\kr_{A}^{\ast}(a_{1})(b))(a_{2})
+\kl_{A^{\ast}}^{\ast}(\kl_{A}^{\ast}(a_{1})(b))(a_{2}).                       \nonumber
\end{align}
\end{pro}

\begin{proof}
It is easy to see that Eq. (\ref{mat1}) is exactly Eq. (\ref{matc1}) and
Eq. (\ref{mat2}) is exactly Eq. (\ref{matc2}) if we take
$\kl_{1}=\kr_{A}^{\ast}-\kl_{A}^{\ast}$, $\kr_{1}=\kr_{A}^{\ast}$,
$\kl_{2}=\kr_{A^{\ast}}^{\ast}-\kl_{A^{\ast}}^{\ast}$, $\kr_{2}=\kr_{A^{\ast}}^{\ast}$.
Thus, Eqs. (\ref{matc1}) and (\ref{matc2}) hold if $(A, A^{\ast},
\kr_{A}^{\ast}-\kl_{A}^{\ast}, \kr_{A}^{\ast}, \kr_{A^{\ast}}^{\ast}-\kl_{A^{\ast}}^{\ast},
\kr_{A^{\ast}}^{\ast})$ is a matched pair of perm algebras. Conversely,
if Eqs. (\ref{matc1}) and (\ref{matc2}) hold, for any $a_{1}, a_{2}\in A$ and
$b_{1}, b_{2}\in A^{\ast}$, since
\begin{align*}
\langle\kr_{A^{\ast}}^{\ast}(b_{1})(a_{1}a_{2}),\; b_{2}\rangle&=\langle a_{1}a_{2},\;
b_{2}b_{1}\rangle=\langle a_{2},\; \kl_{A}^{\ast}(a_{1})(b_{2}b_{1})\rangle,\\
\langle\kl_{A^{\ast}}^{\ast}(b_{1})(a_{1}a_{2}),\; b_{2}\rangle&=\langle a_{1}a_{2},\;
b_{1}b_{2}\rangle=\langle a_{2},\; \kl_{A}^{\ast}(a_{1})(b_{1}b_{2})\rangle,\\
\langle\kr_{A^{\ast}}^{\ast}(b_{1})(a_{2}a_{1}),\; b_{2}\rangle&=\langle a_{2}a_{1},\;
b_{2}b_{1}\rangle=\langle a_{2},\; \kr_{A}^{\ast}(a_{1})(b_{2}b_{1})\rangle,\\
\langle\kl_{A^{\ast}}^{\ast}(b_{1})(a_{2}a_{1}),\; b_{2}\rangle&=\langle a_{2}a_{1},\;
b_{1}b_{2}\rangle=\langle a_{2},\; \kr_{A}^{\ast}(a_{1})(b_{1}b_{2})\rangle.
\end{align*}
we get $\kr_{A}^{\ast}(a)(b_{1}b_{2})-\kl_{A}^{\ast}(a)(b_{1}b_{2})
=\kr_{A}^{\ast}(a)(b_{2}b_{1})-\kl_{A}^{\ast}(a)(b_{2}b_{1})$.
Moreover, by Eq. (\ref{matc2}), we have
\begin{align*}
&\; \kr_{A^{\ast}}^{\ast}(b)(a_{1}a_{2})-\kr_{A^{\ast}}^{\ast}(b)(a_{2}a_{1})\\
=&\; a_{1}\kr_{A^{\ast}}^{\ast}(b)(a_{2})-\kr_{A^{\ast}}^{\ast}(b)(a_{2})a_{1}
+\kl_{A^{\ast}}^{\ast}(\kr_{A}^{\ast}(a_{2})(b))(a_{1})
-\kl_{A^{\ast}}^{\ast}(\kl_{A}^{\ast}(a_{2})(b))(a_{1}).
\end{align*}
Since for any $a_{1}, a_{2}\in A$ and $b_{1}, b_{2}\in A^{\ast}$,
\begin{align*}
\langle\kr_{A^{\ast}}^{\ast}(b_{1})(a_{1}a_{2}),\; b_{2}\rangle&=\langle a_{1}a_{2},\;
b_{2}b_{1}\rangle=\langle a_{2},\; \kl_{A}^{\ast}(a_{1})(b_{2}b_{1})\rangle,\\
\langle\kr_{A^{\ast}}^{\ast}(b_{1})(a_{2}a_{1}),\; b_{2}\rangle&=\langle a_{2}a_{1},\;
b_{2}b_{1}\rangle=\langle a_{2},\; \kr_{A}^{\ast}(a_{1})(b_{2}b_{1})\rangle,\\
\langle a_{1}\kr_{A^{\ast}}^{\ast}(b_{1})(a_{2}),\; b_{2}\rangle
&=\langle\kr_{A^{\ast}}^{\ast}(b_{1})(a_{2}),\; \kl_{A}^{\ast}(a_{1})(b_{2})\rangle
=\langle a_{2},\; \kl_{A}^{\ast}(a_{1})(b_{2})b_{1}\rangle,\\
\langle\kr_{A^{\ast}}^{\ast}(b_{1})(a_{2})a_{1},\; b_{2}\rangle
&=\langle\kr_{A^{\ast}}^{\ast}(b_{1})(a_{2}),\; \kr_{A}^{\ast}(a_{1})(b_{2})\rangle
=\langle a_{2},\; \kr_{A}^{\ast}(a_{1})(b_{2})b_{1}\rangle,\\
\langle\kl_{A^{\ast}}^{\ast}(\kr_{A}^{\ast}(a_{2})(b_{1}))(a_{1}),\; b_{2}\rangle
&=\langle\kr_{A}^{\ast}(a_{2})(b_{1}),\; \kr_{A^{\ast}}^{\ast}(b_{2})(a_{1})\rangle
=\langle a_{2},\; \kl_{A}^{\ast}(\kr_{A^{\ast}}^{\ast}(b_{2})(a_{1}))(b_{1})\rangle,\\
\langle\kl_{A^{\ast}}^{\ast}(\kl_{A}^{\ast}(a_{2})(b_{1}))(a_{1}),\; b_{2}\rangle
&=\langle\kl_{A}^{\ast}(a_{2})(b_{1}),\; \kr_{A^{\ast}}^{\ast}(b_{2})(a_{1})\rangle
=\langle a_{2},\; \kr_{A}^{\ast}(\kr_{A^{\ast}}^{\ast}(b_{2})(a_{1}))(b_{1})\rangle,
\end{align*}
we get
\begin{align*}
&\; \kr_{A}^{\ast}(a)(b_{1}b_{2})-\kl_{A}^{\ast}(a)(b_{1}b_{2})\\
=&\; \kr_{A}^{\ast}(a)(b_{1})b_{2}-\kl_{A}^{\ast}(a)(b_{1})b_{2}
+\kr_{A}^{\ast}(\kr_{A^{\ast}}^{\ast}(b_{1})(a))(b_{2})
-\kl_{A}^{\ast}(\kr_{A^{\ast}}^{\ast}(b_{1})(a))(b_{2}).
\end{align*}
This implies that Eq. (\ref{mat3}) holds. Similarly, one can check
that Eq. (\ref{mat4}) also holds. Hence, $(A, A^{\ast},
\kr_{A}^{\ast}-\kl_{A}^{\ast}, \kr_{A}^{\ast}, \kr_{A^{\ast}}^{\ast}-\kl_{A^{\ast}}^{\ast},
\kr_{A^{\ast}}^{\ast})$ is a matched pair of perm algebras. The proof is finished.
\end{proof}

\begin{lem}\label{lem:cosp}
Let $A$ be a vector space and $\Delta: A\rightarrow A\otimes A$ be a linear map. Then
the dual map $\Delta^{\ast}: A^{\ast}\otimes A^{\ast}\rightarrow A^{\ast}$ defines a
perm algebra structure on $A^{\ast}$ if and only if $\Delta$ satisfies
\begin{align}
(\Delta\otimes\id)\Delta(a)=(\id\otimes\Delta)\Delta(a)
=(\id\otimes\tau\Delta)\Delta(a),                        \label{coalg}
\end{align}
for any $a\in A$, where the twisting map $\tau: A\otimes A\rightarrow A\otimes A$ is
given by $\tau(a_{1}\otimes a_{2})=a_{2}\otimes a_{1}$ for all $a_{1}, a_{2}\in A$.
\end{lem}

\begin{proof}
For any $b_{1}, b_{2}, b_{3}\in A^{\ast}$, since
\begin{align*}
\langle(\Delta\otimes\id)\Delta(a),\; b_{1}\otimes b_{2}\otimes b_{3}\rangle
& =\langle a,\; \Delta^{\ast}(\Delta^{\ast}(b_{1}\otimes b_{2})\otimes b_{3})\rangle, \\
\langle(\id\otimes\Delta)\Delta(a),\; b_{1}\otimes b_{2}\otimes b_{3}\rangle
& =\langle a,\; \Delta^{\ast}(b_{1}\otimes\Delta^{\ast}(b_{2}\otimes b_{3}))\rangle,\\
\langle(\id\otimes\tau\Delta)\Delta(a),\; b_{1}\otimes b_{2}\otimes b_{3}\rangle
& =\langle a,\; \Delta^{\ast}(b_{1}\otimes\Delta^{\ast}(b_{3}\otimes b_{2}))\rangle,
\end{align*}
we get this lemma.
\end{proof}

\begin{defi}\label{def: perm-bi}
Let $A$ be a perm algebra. A {\rm perm bialgebra} structure on $A$ is a linear map
$\Delta: A\rightarrow A\otimes A$ such that $\Delta$ satisfies Eq. (\ref{coalg}), and
\begin{align}
&\Delta(a_{1}a_{2})-\tau\Delta(a_{1}a_{2})
=\Delta(a_{2}a_{1})-\tau\Delta(a_{2}a_{1})                                    \label{bialg1}\\
&\quad =(\kr_{A}(a_{2})\otimes\id)\Delta(a_{1})
-\tau(\id\otimes\kr_{A}(a_{2}))\Delta(a_{1})
+(\id\otimes\kr_{A}(a_{1}))\Delta(a_{2})
-\tau(\kr_{A}(a_{1})\otimes\id)\Delta(a_{2}),                                     \nonumber\\
& \Delta(a_{1}a_{2})=(\kl_{A}(a_{1})\otimes\id)\Delta(a_{2})
+(\id\otimes\kr_{A}(a_{2}))\Delta(a_{1})
-(\id\otimes\kl_{A}(a_{2}))\Delta(a_{1})                                \label{bialg2}\\
&\quad =(\kl_{A}(a_{1})\otimes\id)\Delta(a_{2})
-\tau(\id\otimes\kl_{A}(a_{1}))\Delta(a_{2})
+(\id\otimes\kr_{A}(a_{2}))\Delta(a_{1})                                       \nonumber\\
&\quad =(\kr_{A}(a_{2})\otimes\id)\Delta(a_{1})
+(\id\otimes\kr_{A}(a_{1}))\Delta(a_{2})
-(\id\otimes\kl_{A}(a_{1}))\Delta(a_{2})                                      \nonumber\\
&\qquad\qquad\qquad\qquad\qquad -\tau(\kr_{A}(a_{1})\otimes\id)\Delta(a_{2})
+\tau(\kl_{A}(a_{1})\otimes\id)\Delta(a_{2}).                                \nonumber
\end{align}
We denote this perm bialgebra by $(A, \Delta)$.
\end{defi}

\begin{ex}\label{ex:bialg}
Let $A=\mathbb{C}\{e_{1}, e_{2}\}$ be a 2-dimensional perm algebra with
non-zero products: $e_{1}e_{1}=e_{1}$ and $e_{2}e_{1}=e_{2}$.
By direct computations, we can obtain that a linear map $\Delta: A\rightarrow A\otimes A$
such that $(A, \Delta)$ is a perm bialgebra if and only if $\Delta(e_{1})=k_{1}
e_{2}\otimes e_{2}$ and $\Delta(e_{2})=k_{2}e_{2}\otimes e_{2}$, for some
$k_{1}, k_{2}\in\mathbb{C}$.
\end{ex}

It is easy to see that Eqs. (\ref{bialg1}) and (\ref{bialg2}) hold if and only if
Eqs. (\ref{matc1}) and (\ref{matc2}) hold. Combining Propositions \ref{pro:main-mat}
and \label{pro:mat-bia}, we have the following conclusion.

\begin{thm}\label{thm: equat}
Let $A$ be a perm algebra and $\Delta: A\rightarrow A\otimes A$ be a linear map such
that $\Delta^{\ast}: A^{\ast}\otimes A^{\ast}\rightarrow A^{\ast}$ defines a perm algebra
structure on $A^{\ast}$.  Then the following conditions are equivalent:
\begin{itemize}
\item[$(i)$] $(A, \Delta)$ is a perm bialgebra;
\item[$(ii)$] $(A, A^{\ast}, \kr_{A}^{\ast}-\kl_{A}^{\ast}, \kr_{A}^{\ast},
     \kr_{A^{\ast}}^{\ast}-\kl_{A^{\ast}}^{\ast}, \kr_{A^{\ast}}^{\ast})$ is a matched
     pair of perm algebras;
\item[$(iii)$] $(A\oplus A^{\ast}, A, A^{\ast})$ is a standard Manin triple with respect
     to the bilinear form $\widehat{\omega}$.
\end{itemize}
\end{thm}

\subsection{Coboundary perm bialgebras and $\mathcal{S}$-equation in perm algebas}
\label{subsec:cobou}

In this subsection, we consider a special class of perm bialgebras, that is, a perm
bialgebra $(A, \Delta)$ with $\Delta$ in the form
\begin{align}
\Delta(a)=(\kl_{A}(a)\otimes\id+\id\otimes\kl_{A}(a)
-\id\otimes\kr_{A}(a))\mathfrak{R},                    \label{rmax}
\end{align}
for any $a\in A$, where $\mathfrak{R}\in A\otimes A$.

Such perm bialgebras are quite similar as the ``coboundary left-symmetric bialgebras"
for left-symmetric algebras \cite{Bai}. We call it {\it coboundary perm bialgebras}.
Moreover, it is easy to see that $\tau\Delta(a)=(\id\otimes\kr_{A}(a)+\kr_{A}(a)
\otimes\id)\tau\mathfrak{R}$ if Eq. (\ref{rmax}) holds.

\begin{pro}\label{pro:cob}
Let $A$ be a perm algebra, $\mathfrak{R}\in A\otimes A$, and $\Delta: A\rightarrow
A\otimes A$ be a linear map defined by Eq. (\ref{rmax}). Then for any $a_{1}, a_{2}\in A$,
we have
\begin{enumerate}
\item[$(i)$] $\Delta$ satisfies Eq. (\ref{bialg1}) if and only if
  \begin{align}
  &\Big(\kl_{A}(a_{1}a_{2})\otimes\id+\id\otimes\kl_{A}(a_{1}a_{2})\Big)
  (\mathfrak{R}-\tau\mathfrak{R})=\Big(\kl_{A}(a_{2}a_{1})\otimes\id
  +\id\otimes\kl_{A}(a_{2}a_{1})\Big)(\mathfrak{R}-\tau\mathfrak{R}),    \label{cob1}\\
  &\qquad\quad\; \Big(\kr_{A}(a_{2})\otimes\kl_{A}(a_{1})
  +\kl_{A}(a_{2})\otimes\kr_{A}(a_{1})+\id\otimes\kl_{A}(a_{2}a_{1})\Big)
  (\mathfrak{R}-\tau\mathfrak{R})                                        \label{cob2}\\[-1mm]
  &\qquad = \Big(\kr_{A}(a_{2})\otimes\kr_{A}(a_{1})
  +\id\otimes\kr_{A}(a_{2}a_{1})+\id\otimes\kl_{A}(a_{1}a_{2})\Big)
  (\mathfrak{R}-\tau\mathfrak{R});                                       \nonumber
  \end{align}
\item[$(ii)$] $\Delta$ satisfies Eq. (\ref{bialg2}) if and only if
  \begin{align}
  &\qquad\qquad\qquad (\kl_{A}(a_{1})\otimes\kl_{A}(a_{2}))
  (\mathfrak{R}-\tau\mathfrak{R})=0,                                    \label{cob3}\\
  &\quad\; \Big(\kr_{A}(a_{2})\otimes\kl_{A}(a_{1})
  +\kl_{A}(a_{2})\otimes\kr_{A}(a_{1})+\id\otimes\kl_{A}(a_{2}a_{1})\Big)
  (\mathfrak{R}-\tau\mathfrak{R})                                        \label{cob4}\\[-1mm]
  & = \Big(\kr_{A}(a_{2})\otimes\kr_{A}(a_{1})
  +\kl_{A}(a_{2})\otimes\kl_{A}(a_{1})+\id\otimes\kl_{A}(a_{1}a_{2})\Big)
  (\mathfrak{R}-\tau\mathfrak{R}).                                       \nonumber
  \end{align}
\end{enumerate}
\end{pro}

\begin{proof}
For any $a_{1}, a_{2}\in A$, we have
\begin{align*}
&\; \Delta(a_{1}a_{2})-\tau\Delta(a_{1}a_{2})-\Delta(a_{2}a_{1})+\tau\Delta(a_{2}a_{1})\\
=&\; \Big(\kl_{A}(a_{1}a_{2})\otimes\id+\id\otimes\kl_{A}(a_{1}a_{2})
-\kl_{A}(a_{2}a_{1})\otimes\id+\id\otimes\kl_{A}(a_{2}a_{1})\Big)
(\mathfrak{R}-\tau\mathfrak{R}).
\end{align*}
Moreover, note that
\begin{align*}
(\kr_{A}(a_{2})\otimes\id)\Delta(a_{1})
& =(\kr_{A}(a_{2})\kl_{A}(a_{1})\otimes\id)\mathfrak{R}
+(\kr_{A}(a_{2})\otimes\kl_{A}(a_{1}))\mathfrak{R}
-(\kr_{A}(a_{2})\otimes\kr_{A}(a_{1}))\mathfrak{R},\\
\tau(\id\otimes\kr_{A}(a_{2}))\Delta(a_{1})
& =(\kr_{A}(a_{2})\otimes\kl_{A}(a_{1}))\tau\mathfrak{R}
+(\kr_{A}(a_{2})\kl_{A}(a_{1})\otimes\id)\tau\mathfrak{R}
-(\kr_{A}(a_{2})\kr_{A}(a_{1})\otimes\id)\tau\mathfrak{R},\\
(\id\otimes\kr_{A}(a_{1}))\Delta(a_{2})
&=(\kl_{A}(a_{2})\otimes\kr_{A}(a_{1}))\mathfrak{R}
+(\id\otimes\kr_{A}(a_{1})\kl_{A}(a_{2}))\mathfrak{R}
-(\id\otimes\kr_{A}(a_{1})\kr_{A}(a_{2}))\mathfrak{R},\\
\tau(\kr_{A}(a_{1})\otimes\id)\Delta(a_{2})
& =(\id\otimes\kr_{A}(a_{1})\kl_{A}(a_{2}))\tau\mathfrak{R}
+(\kl_{A}(a_{2})\otimes\kr_{A}(a_{1}))\tau\mathfrak{R}
-(\kr_{A}(a_{2})\otimes\kr_{A}(a_{1}))\tau\mathfrak{R},
\end{align*}
we obtain
\begin{align*}
&\; (\kr_{A}(a_{2})\otimes\id)\Delta(a_{1})-\tau(\id\otimes\kr_{A}(a_{2}))\Delta(a_{1})
+(\id\otimes\kr_{A}(a_{1}))\Delta(a_{2})\\
&\qquad\qquad\qquad\qquad -\tau(\kr_{A}(a_{1})\otimes\id)\Delta(a_{2})
+\tau\Delta(a_{1}a_{2})-\Delta(a_{1}a_{2})\\
=&\; \Big(\kr_{A}(a_{2})\otimes\kl_{A}(a_{1})+\kl_{A}(a_{1})\otimes\kr_{A}(a_{1})
+\id\otimes\kl_{A}(a_{2}a_{1})-\kr_{A}(a_{2})\otimes\kr_{A}(a_{1})\\[-1mm]
&\qquad\qquad\qquad\qquad\qquad\qquad-\id\otimes\kr_{A}(a_{2}a_{1})
-\id\otimes\kl_{A}(a_{1}a_{2})\Big)(\mathfrak{R}-\tau\mathfrak{R}).
\end{align*}
Therefore, Eq. (\ref{bialg1}) holds if and only if Eqs. (\ref{cob1}) and (\ref{cob2}) hold.
Similarly, Eq. (\ref{bialg2}) holds if and only if Eqs. (\ref{cob3})-(\ref{cob4}) hold.
\end{proof}

Let $A$ be a perm algebra and $\mathfrak{R}=\sum_{i}a_{i}\otimes\bar{a}_{i}\in A\otimes A$.
Set
\begin{align*}
\mathfrak{R}_{12}&=\sum_{i}a_{i}\otimes\bar{a}_{i}\otimes{\bf1},\qquad
\mathfrak{R}_{13}=\sum_{i}a_{i}\otimes{\bf1}\otimes\bar{a}_{i},\qquad
\mathfrak{R}_{23}=\sum_{i}{\bf1}\otimes a_{i}\otimes\bar{a}_{i},\\[-1mm]
\mathfrak{R}_{21}&=\sum_{i}\bar{a}_{i}\otimes a_{i}\otimes{\bf1},\qquad
\mathfrak{R}_{31}=\sum_{i}\bar{a}_{i}\otimes{\bf1}\otimes a_{i},\qquad
\mathfrak{R}_{32}=\sum_{i}{\bf1}\otimes\bar{a}_{i}\otimes a_{i},
\end{align*}
where ${\bf1}$ is the unit if $A$ has a unit, otherwise is a symbol playing a similar
role of the unit. Then the operation between two $\mathfrak{R}$s is in an obvious way.
For example,
$$
\mathfrak{R}_{12}\mathfrak{R}_{13}=\sum_{i,j}a_{i}a_{j}\otimes
\bar{a}_{i}\otimes\bar{a}_{j},\quad
\mathfrak{R}_{13}\mathfrak{R}_{23}=\sum_{i,j}a_{i}\otimes a_{j}
\otimes\bar{a}_{i}\bar{a}_{j},\quad
\mathfrak{R}_{23}\mathfrak{R}_{12}=\sum_{i,j}a_{j}\otimes a_{i}
\bar{a}_{j}\otimes\bar{a}_{i},
$$
and so on.

\begin{pro}\label{pro:dual-coalg}
Let $A$ be a perm algebra, $\mathfrak{R}=\sum_{i}a_{i}\otimes\bar{a}_{i}\in A\otimes A$,
and $\Delta: A\rightarrow A\otimes A$ be a linear map defined by Eq. (\ref{rmax}).
Then $\Delta^{\ast}$ defines a perm algebra structure on $A^{\ast}$ if and only if
for any $a\in A$,
\begin{align}
(\id\otimes\id\otimes(\kl_{A}(a)-\kr_{A}(a)))P(\mathfrak{R})
&=(\kl_{A}(a)\otimes\id\otimes\id)S(\mathfrak{R}),                \label{coalg1}\\
(\id\otimes\id\otimes(\kl_{A}(a)-\kr_{A}(a)))Q(\mathfrak{R})
&=(\kl_{A}(a)\otimes\id\otimes\id)T(\mathfrak{R})+M(\mathfrak{R}), \label{coalg2}
\end{align}
where
\begin{align*}
P(\mathfrak{R})&=\mathfrak{R}_{13}\mathfrak{R}_{12}-\mathfrak{R}_{13}\mathfrak{R}_{23}
+[\mathfrak{R}_{23}, \mathfrak{R}_{12}],\\
Q(\mathfrak{R})&=\mathfrak{R}_{13}\mathfrak{R}_{12}-\mathfrak{R}_{13}\mathfrak{R}_{32}
+[\mathfrak{R}_{23}, \mathfrak{R}_{12}],\\
S(\mathfrak{R})&=\mathfrak{R}_{12}\mathfrak{R}_{23}+\mathfrak{R}_{13}\mathfrak{R}_{23}
-\mathfrak{R}_{12}\mathfrak{R}_{13}-\mathfrak{R}_{23}\mathfrak{R}_{13},\\
T(\mathfrak{R})&=\mathfrak{R}_{12}\mathfrak{R}_{23}+\mathfrak{R}_{13}\mathfrak{R}_{32}
-\mathfrak{R}_{13}\mathfrak{R}_{12}-\mathfrak{R}_{32}\mathfrak{R}_{12},\\
M(\mathfrak{R})&=(\id\otimes\kl_{A}(a)\otimes\id)(\mathfrak{R}_{32}\mathfrak{R}_{12}
-\mathfrak{R}_{23}\mathfrak{R}_{12})+(\id\otimes\kr_{A}(a)\otimes\id)
(\mathfrak{R}_{12}\mathfrak{R}_{23}-\mathfrak{R}_{12}\mathfrak{R}_{32}),
\end{align*}
and $[\mathfrak{R}_{23}, \mathfrak{R}_{12}]=\mathfrak{R}_{23}\mathfrak{R}_{12}
-\mathfrak{R}_{12}\mathfrak{R}_{23}$.
\end{pro}

\begin{proof}
For any $a\in A$, we have
\begin{align*}
(\Delta\otimes\id)\Delta(a)
&=\sum_{i,j}\Big((aa_{i})a_{j}\otimes\bar{a}_{j}\otimes\bar{a}_{i}
+a_{j}\otimes(aa_{i})\bar{a}_{j}\otimes\bar{a}_{i}
-a_{j}\otimes\bar{a}_{j}(aa_{i})\otimes\bar{a}_{i}\\[-4mm]
&\qquad\quad+a_{i}a_{j}\otimes\bar{a}_{j}\otimes a\bar{a}_{i}
+a_{j}\otimes a_{i}\bar{a}_{j}\otimes a\bar{a}_{i}
-a_{j}\otimes \bar{a}_{j}a_{i}\otimes a\bar{a}_{i}\\[-2mm]
&\qquad\qquad-a_{i}a_{j}\otimes\bar{a}_{j}\otimes\bar{a}_{i}a
-a_{j}\otimes a_{i}\bar{a}_{j}\otimes\bar{a}_{i}a
+a_{j}\otimes \bar{a}_{j}a_{i}\otimes\bar{a}_{i}a\Big),\\
%
(\id\otimes\Delta)\Delta(a)
&=\sum_{i,j}\Big(aa_{i}\otimes\bar{a}_{i}a_{j}\otimes\bar{a}_{j}
+aa_{i}\otimes a_{j}\otimes\bar{a}_{i}\bar{a}_{j}
-aa_{i}\otimes a_{j}\otimes\bar{a}_{j}\bar{a}_{i}\\[-4mm]
&\qquad\quad +a_{i}\otimes(a\bar{a}_{i})a_{j}\otimes\bar{a}_{j}
+a_{i}\otimes a_{j}\otimes(a\bar{a}_{i})\bar{a}_{j}
-a_{i}\otimes a_{j}\otimes\bar{a}_{j}(a\bar{a}_{i})\\[-2mm]
&\qquad\qquad -a_{i}\otimes(\bar{a}_{i}a)a_{j}\otimes\bar{a}_{j}
-a_{i}\otimes a_{j}\otimes(\bar{a}_{i}a)\bar{a}_{j}
+a_{i}\otimes a_{j}\otimes\bar{a}_{j}(\bar{a}_{i}a)\Big),\\
%
(\id\otimes\tau\Delta)\Delta(a)
&=\sum_{i,j}\Big(aa_{i}\otimes\bar{a}_{j}\otimes\bar{a}_{i}a_{j}
+aa_{i}\otimes\bar{a}_{i}\bar{a}_{j}\otimes a_{j}
-aa_{i}\otimes\bar{a}_{j}\bar{a}_{i}\otimes a_{j}\\[-4mm]
&\qquad\quad +a_{i}\otimes\bar{a}_{j}\otimes(a\bar{a}_{i})a_{j}
+a_{i}\otimes(a\bar{a}_{i})\bar{a}_{j}\otimes a_{j}
-a_{i}\otimes\bar{a}_{j}(a\bar{a}_{i})\otimes a_{j}\\[-2mm]
&\qquad\qquad -a_{i}\otimes\bar{a}_{j}\otimes(\bar{a}_{i}a)a_{j}
-a_{i}\otimes(\bar{a}_{i}a)\bar{a}_{j}\otimes a_{j}
+a_{i}\otimes\bar{a}_{j}(\bar{a}_{i}a)\otimes a_{j}\Big).
\end{align*}
Thus, $(\Delta\otimes\id)\Delta(a)=(\id\otimes\Delta)\Delta(a)$ is equivalent to
\begin{align*}
&\;(\id\otimes\id\otimes(\kl_{A}(a)-\kr_{A}(a)))\big(\mathfrak{R}_{13}\mathfrak{R}_{12}
-\mathfrak{R}_{13}\mathfrak{R}_{23}+[\mathfrak{R}_{23}, \mathfrak{R}_{12}]\big)\\
=&\;(\kl_{A}(a)\otimes\id\otimes\id)(\mathfrak{R}_{12}\mathfrak{R}_{23}
+\mathfrak{R}_{13}\mathfrak{R}_{23}-\mathfrak{R}_{12}\mathfrak{R}_{13}
-\mathfrak{R}_{23}\mathfrak{R}_{13}),
\end{align*}
and $(\id\otimes\Delta)\Delta(a)=(\id\otimes\tau\Delta)\Delta(a)$ is equivalent to
\begin{align*}
&\;(\id\otimes\id\otimes(\kl_{A}(a)-\kr_{A}(a)))\big(\mathfrak{R}_{13}\mathfrak{R}_{12}
-\mathfrak{R}_{13}\mathfrak{R}_{32}+[\mathfrak{R}_{23}, \mathfrak{R}_{12}]\big)\\
=&\; (\kl_{A}(a)\otimes\id\otimes\id)(\mathfrak{R}_{12}\mathfrak{R}_{23}
+\mathfrak{R}_{13}\mathfrak{R}_{32}-\mathfrak{R}_{13}\mathfrak{R}_{12}
-\mathfrak{R}_{32}\mathfrak{R}_{12})\\
&\quad+(\id\otimes\kl_{A}(a)\otimes\id)(\mathfrak{R}_{32}\mathfrak{R}_{12}
-\mathfrak{R}_{23}\mathfrak{R}_{12})+(\id\otimes\kr_{A}(a)\otimes\id)
(\mathfrak{R}_{12}\mathfrak{R}_{23}-\mathfrak{R}_{12}\mathfrak{R}_{32}).
\end{align*}
The proof is finished.
\end{proof}

Combing Propositions \ref{pro:cob} and \ref{pro:dual-coalg} together,
we have the following conclusion.

\begin{thm}\label{thm:bialg-perm}
Let $A$ be a perm algebra, $\mathfrak{R}\in A\otimes A$,
and $\Delta: A\rightarrow A\otimes A$ be a linear map defined by Eq. (\ref{rmax}).
Then $(A, \Delta)$ is a perm bialgebra if and only if $\mathfrak{R}$ satisfies Eqs.
(\ref{cob1})-(\ref{coalg2}).
\end{thm}

For any $a_{1}, a_{2}, a_{3}\in A$, we define $\tau_{23}(a_{1}\otimes a_{2}\otimes a_{3})
=a_{1}\otimes a_{3}\otimes a_{2}$. Then it is easy to see that $S(\mathfrak{R})=
\tau_{23}T(\mathfrak{R})$. If $\mathfrak{R}\in A\otimes A$ is {\it symmetric}, i.e.,
$\mathfrak{R}=\tau\mathfrak{R}$, then the Eqs. (\ref{cob1})-(\ref{cob4}) hold,
$P(\mathfrak{R})=Q(\mathfrak{R})=-T(\mathfrak{R})$ and $M(\mathfrak{R})=0$. Thus
we have the following corollary.

\begin{cor}\label{cor:bialg-perm}
Let $A$ be a perm algebra, $\mathfrak{R}\in A\otimes A$ be symmetric,
and $\Delta: A\rightarrow A\otimes A$ be a linear map defined by Eq. (\ref{rmax}). If
\begin{align}
P(\mathfrak{R})=\mathfrak{R}_{13}\mathfrak{R}_{12}-\mathfrak{R}_{13}\mathfrak{R}_{23}
+[\mathfrak{R}_{23}, \mathfrak{R}_{12}]=0,                     \label{YBE}
\end{align}
then $(A, \Delta)$ is a perm bialgebra.
\end{cor}

\begin{defi}\label{def:YBE}
Let $A$ be a perm algebra and $\mathfrak{R}\in A\otimes A$. The Eq. (\ref{YBE}) is called the
{\rm $\mathcal{S}$-equation} in $A$.
\end{defi}

The notion of $\mathcal{S}$-equation in a perm algebra, similar to the
$\mathcal{S}$-equation in a left-symmetric algebra, which is due to the fact that
it is an analogue of the classical Yang-Baxter equation in a Lie algebra or the
associative Yang-Baxter equation in an associative algebra.

\begin{ex}\label{ex: s-equ}
Let $A=\mathbb{C}\{e_{1}, e_{2}\}$ be a 2-dimensional perm algebra with
non-zero products: $e_{1}e_{1}=e_{1}$ and $e_{2}e_{1}=e_{2}$.
Consider the symmetric element $\mathfrak{R}=e_{1}\otimes e_{2}+e_{2}\otimes e_{1}
+e_{2}\otimes e_{2}\in A\otimes A$. One can check that $\mathfrak{R}$ satisfies Eqs.
(\ref{cob1})-(\ref{coalg2}). Thus the Eq. (\ref{rmax}) give a perm bialgebra structure
on $A$ by $\Delta(e_{1})=e_{2}\otimes e_{2}$ and $\Delta(e_{2})=2e_{2}\otimes e_{2}$.
But in this case $P(\mathfrak{R})=\mathfrak{R}_{13}\mathfrak{R}_{12}
-\mathfrak{R}_{13}\mathfrak{R}_{23}+\mathfrak{R}_{23}\mathfrak{R}_{12}
-\mathfrak{R}_{12}\mathfrak{R}_{23}=e_{2}\otimes e_{2}\otimes e_{1}
-e_{2}\otimes e_{1}\otimes e_{2}\neq0$. This means that $\mathfrak{R}$ is not a
solution of $\mathcal{S}$-equation in $A$. On the other hand, let $\mathfrak{R}=
e_{2}\otimes e_{2}$. It is easy to see that $\mathfrak{R}$ is a symmetric solution
of $\mathcal{S}$-equation in $A$. Therefore, the Eq. (\ref{rmax}) give a perm bialgebra
structure on $A$ by $\Delta(e_{1})=e_{2}\otimes e_{2}$ and $\Delta(e_{2})=0$.
\end{ex}

Next we give some properties of the $\mathcal{S}$-equation in perm algebas.
Let $V$ be a $\bk$-vector space. For all $\mathfrak{R}\in V\otimes V$, we define
linear map $\mathfrak{R}^{\sharp}: V^{\ast}\rightarrow V$ by
$$
\langle\mathfrak{R}^{\sharp}(u_{1}),\; u_{2}\rangle
=\langle\mathfrak{R},\; u_{1}\otimes u_{2}\rangle,
$$
for any $u_{1}\otimes u_{2}\in V^{\ast}$. We say that $\mathfrak{R}\in V\otimes V$ is
{\it nondegenerate}, if the linear map $\mathfrak{R}^{\sharp}$ is an isomorphism.

\begin{pro}  \label{pro:biline}
Let $A$ be a perm algebra and $\mathfrak{R}\in A\otimes A$ be symmetric.
If $\mathfrak{R}$ is nondegenerate, then $\mathfrak{R}$ is a solution of
$\mathcal{S}$-equation in $A$ if and only if the inverse of the isomorphism
$\mathfrak{R}^{\sharp}: A^{\ast}\rightarrow A$, regarded as a bilinear form
$\bar{\omega}$ on $A$, that is, $\bar{\omega}(a_{1}, a_{2})=\langle a_{1},\;
(\mathfrak{R}^{\sharp})^{-1}(a_{2})\rangle$, for any $a_{1}, a_{2}\in A$, satisfies
\begin{align*}
\bar{\omega}(a_{1}a_{2}, a_{3})+\bar{\omega}(a_{1}a_{3}, a_{2})
=\bar{\omega}(a_{3}a_{2}, a_{1})+\bar{\omega}(a_{3}a_{1}, a_{2}),
\end{align*}
for all $a_{1}, a_{2}, a_{3}\in A$.
\end{pro}

\begin{proof}
Let $\mathfrak{R}=\sum_{i}a_{i}\otimes\bar{a}_{i}\in A\otimes A$ be symmetric.
Then $\mathfrak{R}^{\sharp}(b)=\sum_{i}\langle b, a_{i}\rangle\bar{a}_{i}
=\sum_{i}\langle b,\bar{a}_{i}\rangle a_{i}$ for any $b\in A^{\ast}$,
and one can check that $\bar{\omega}$ is a symmetric bilinear form on $A$.
Since $\mathfrak{R}$ is nondegenerate, for any $a_{1}, a_{2}, a_{3}\in A$,
there exist $b_{1}, b_{2}, b_{3}\in A^{\ast}$ such that $a_{j}=\mathfrak{R}^{\sharp}(b_{j})$,
$j=1, 2, 3$. Thus, we have
\begin{align*}
\bar{\omega}(a_{1}a_{2}, a_{3})
&=\langle\mathfrak{R}^{\sharp}(b_{1})\mathfrak{R}^{\sharp}(b_{2}),\; b_{3}\rangle
=\sum_{i, j}\langle b_{1},\; a_{i}\rangle\langle b_{2},\; a_{j}\rangle
\langle b_{3},\; \bar{a}_{i}\bar{a}_{j}\rangle\\[-2mm]
&=\langle b_{1}\otimes b_{2}\otimes b_{3},\; \mathfrak{R}_{13}\mathfrak{R}_{23}\rangle,\\
\bar{\omega}(a_{1}a_{3}, a_{2})
&=\langle\mathfrak{R}^{\sharp}(b_{1})\mathfrak{R}^{\sharp}(b_{3}),\; b_{2}\rangle
=\sum_{i, j}\langle b_{1},\; a_{i}\rangle\langle b_{3},\; a_{j}\rangle
\langle b_{2},\; \bar{a}_{i}\bar{a}_{j}\rangle\\[-2mm]
&=\langle b_{1}\otimes b_{2}\otimes b_{3},\; \mathfrak{R}_{12}\mathfrak{R}_{23}\rangle,\\
\bar{\omega}(a_{3}a_{1}, a_{2})
&=\langle\mathfrak{R}^{\sharp}(b_{3})\mathfrak{R}^{\sharp}(b_{1}),\; b_{2}\rangle
=\sum_{i, j}\langle b_{3},\; a_{i}\rangle\langle b_{1},\; a_{j}\rangle
\langle b_{2},\; \bar{a}_{i}\bar{a}_{j}\rangle\\[-2mm]
&=\langle b_{1}\otimes b_{2}\otimes b_{3},\; \mathfrak{R}_{23}\mathfrak{R}_{12}\rangle,\\
\bar{\omega}(a_{3}a_{2}, a_{1})
&=\langle\mathfrak{R}^{\sharp}(b_{3})\mathfrak{R}^{\sharp}(b_{2}),\; b_{1}\rangle
=\sum_{i, j}\langle b_{3},\; a_{i}\rangle\langle b_{2},\; a_{j}\rangle
\langle b_{1},\; \bar{a}_{i}\bar{a}_{j}\rangle\\[-2mm]
&=\langle b_{1}\otimes b_{2}\otimes b_{3},\; \mathfrak{R}_{13}\mathfrak{R}_{12}\rangle.
\end{align*}
Therefore, $P(\mathfrak{R})=0$ if and only if $\langle b_{1}\otimes b_{2}\otimes b_{3},\;
P(\mathfrak{R})\rangle=0$ for all $b_{1}, b_{2}, b_{3}\in A^{\ast}$, if and only if
$\bar{\omega}(a_{1}a_{2}, a_{3})+\bar{\omega}(a_{1}a_{3}, a_{2})
=\bar{\omega}(a_{3}a_{2}, a_{1})+\bar{\omega}(a_{3}a_{1}, a_{2})$ for all
$a_{1}, a_{2}, a_{3}\in A$.
\end{proof}

For more general symmetric solution of $\mathcal{S}$-equation in a perm algebra, we have the
following Proposition.

\begin{pro}\label{pro:solut}
Let $A$ be a perm algebra and $\mathfrak{R}\in A\otimes A$ be symmetric.
Then $\mathfrak{R}$ is solution of $\mathcal{S}$-equation in $A$ if and only if
\begin{align*}
\mathfrak{R}^{\sharp}(b_{1})\mathfrak{R}^{\sharp}(b_{2})
=\mathfrak{R}^{\sharp}\Big(\kr_{A}^{\ast}(\mathfrak{R}^{\sharp}(b_{1}))(b_{2})
-\kl_{A}^{\ast}(\mathfrak{R}^{\sharp}(b_{1}))(b_{2})
+\kr_{A}^{\ast}(\mathfrak{R}^{\sharp}(b_{2}))(b_{1})\Big),
\end{align*}
for any $b_{1}, b_{2}\in A^{\ast}$.
\end{pro}

\begin{proof}
Let $\{e_{1}, e_{2}, \cdots, e_{n}\}$ be a basis of $A$ and $\{e^{\ast}_{1}, e^{\ast}_{2},
\cdots, e^{\ast}_{n}\}$ be its dual basis. Suppose $e_{i}e_{j}=\sum_{k}
c_{i,j}^{k}e_{k}$ for $1\leq i, j\leq n$, and $\mathfrak{R}=\sum_{i, j}d_{i,j}e_{i}
\otimes e_{j}$, where $d_{i,j}=d_{j,i}$. Then $\mathfrak{R}^{\sharp}(e^{\ast}_{l})
=\sum_{t}d_{tl}e_{t}$. Note that $\mathfrak{R}$ is a solution of $\mathcal{S}$-equation
in $A$ if and only if, for any $1\leq i, j, m\leq n$,
\begin{align*}
\sum_{k, l}\Big(c_{k,l}^{m}d_{i,k}d_{j,l}-c_{k,l}^{j}d_{i,l}d_{k,m}
+c_{l,k}^{i}d_{j,l}d_{k,m}-c_{k,l}^{i}d_{j,l}d_{k,m}\Big)=0.
\end{align*}
The left-hand side of the above equation is just the coefficient of $e_{m}$ in
$$
\mathfrak{R}^{\sharp}(e^{\ast}_{i})\mathfrak{R}^{\sharp}(e^{\ast}_{j})
-\mathfrak{R}^{\sharp}\Big(\kr_{A}^{\ast}(\mathfrak{R}^{\sharp}(e^{\ast}_{i}))(e^{\ast}_{j})
-\kl_{A}^{\ast}(\mathfrak{R}^{\sharp}(e^{\ast}_{i}))(e^{\ast}_{j})
+\kr_{A}^{\ast}(\mathfrak{R}^{\sharp}(e^{\ast}_{j}))(e^{\ast}_{i})\Big).
$$
Since $\mathfrak{R}^{\sharp}$, $\kl_{A}^{\ast}$ and $\kr_{A}^{\ast}$ are linear,
the conclusion follows.
\end{proof}

\begin{pro}\label{pro: iso}
Let $A$ be a perm algebra, $\mathfrak{R}\in A\otimes A$ be a symmetric
solution of $\mathcal{S}$-equation in $A$, and $\Delta: A\rightarrow A\otimes A$
be a linear map defined by Eq. (\ref{rmax}).
\begin{enumerate}
\item[$(i)$] the perm algebra structure induced by $\Delta$ on the dual space $A^{\ast}$ is
  given by
  $$
  b_{1}b_{2}=\kr_{A}^{\ast}(\mathfrak{R}^{\sharp}(b_{2}))(b_{1})
  +\kr_{A}^{\ast}(\mathfrak{R}^{\sharp}(b_{1}))(b_{2})
  -\kl_{A}^{\ast}(\mathfrak{R}^{\sharp}(b_{1}))(b_{2}),
  $$
  for any $b_{1}, b_{2}\in A^{\ast}$;

\item[$(ii)$] if the perm algebra structure on $A^{\ast}$ is induced by $\Delta$,
  then $\mathfrak{R}^{\sharp}: A^{\ast}\rightarrow A$ is a morphism of perm algebras;
  if the addition of $\mathfrak{R}$ is nondegenerate, then $\mathfrak{R}^{\sharp}$
  is a perm algebra isomorphism.
\end{enumerate}
\end{pro}

\begin{proof}
Let $\{e_{1}, e_{2}, \cdots, e_{n}\}$ be a basis of $A$ and $\{e^{\ast}_{1}, e^{\ast}_{2},
\cdots, e^{\ast}_{n}\}$ be its dual basis. Suppose $e_{i}e_{j}=\sum_{k}
c_{i,j}^{k}e_{k}$, $e^{\ast}_{i}e^{\ast}_{j}=\sum_{k}f_{i,j}^{k}e^{\ast}_{k}$ for
$1\leq i, j\leq n$, and $\mathfrak{R}=\sum_{i, j}d_{i,j}e_{i}\otimes e_{j}$, where
$d_{i,j}=d_{j,i}$. Then $\mathfrak{R}^{\sharp}(e^{\ast}_{l})=\sum_{t}d_{tl}e_{t}$.
Moreover, since
$$
\Delta(e_{i})=\sum_{k,l}f^{i}_{k,l}e_{k}\otimes e_{l}
=(\kl_{A}(e_{i})\otimes\id+\id\otimes\kl_{A}(e_{i})-\id\otimes\kr_{A}(e_{i}))\mathfrak{R},
$$
we have $f^{i}_{k,l}=\sum_{t}\big(d_{t,l}c_{i,t}^{k}+d_{k,t}c_{i,t}^{l}
-d_{k,t}c_{t,i}^{l}\big)$ for all $1\leq i, j, k\leq n$. Thus,
\begin{align*}
e^{\ast}_{k}e^{\ast}_{l}
&=\sum_{i, t}\Big(d_{t,l}c_{i,t}^{k}+d_{k,t}c_{i,t}^{l}
-d_{k,t}c_{t,i}^{l}\Big)e^{\ast}_{i}\\[-2mm]
&=\sum_{i, t}\Big(d_{t,l}\langle e^{\ast}_{k},\; e_{i}e_{t}\rangle
+d_{k,t}\langle e^{\ast}_{l},\; e_{i}e_{t}\rangle
-d_{k,t}\langle e^{\ast}_{l},\; e_{t}e_{i}\rangle\Big)e^{\ast}_{i}\\[-2mm]
&=\sum_{i}\Big(\langle e^{\ast}_{k},\; e_{i}\mathfrak{R}^{\sharp}(e^{\ast}_{l})\rangle
+\langle e^{\ast}_{l},\; e_{i}\mathfrak{R}^{\sharp}(e^{\ast}_{k})\rangle
-\langle e^{\ast}_{l},\; \mathfrak{R}^{\sharp}(e^{\ast}_{k})
e_{i}\rangle\Big)e^{\ast}_{i}\\[-2mm]
&=\kr_{A}^{\ast}(\mathfrak{R}^{\sharp}(e^{\ast}_{l}))(e^{\ast}_{k})
+\kr_{A}^{\ast}(\mathfrak{R}^{\sharp}(e^{\ast}_{k}))(e^{\ast}_{l})
-\kl_{A}^{\ast}(\mathfrak{R}^{\sharp}(e^{\ast}_{k}))(e^{\ast}_{l}).
\end{align*}
Thus, we get $(i)$. The conclusion $(ii)$ following from $(i)$ and Proposition \ref{pro:solut}.
\end{proof}

\begin{rmk}
In this paper, we study the extending structures problem for perm algebras,
perm bialgebras, and some issues related to perm algebras. In the end of this paper,
we collect two further questions for the perm algebas. First, the cohomology theory.
The perm operad is Koszul. How to clearly define the cohomology
complex of perme algebras and the differential graded Lie algebra structure on cohomology
complex to study the deformation theory of perm algebras?
Second, the $\mathcal{O}$-operator on perm algebras. The $\mathcal{O}$-operator
on Lie algebra is a natural generalization of the classical Yang-Baxter equation in
a Lie algebra. The former gives a construction of skew-symmetric solutions of classical
Yang-Baxter equation. For the $\mathcal{O}$-operator on a perm algebra and related problem,
we will investigate in following paper \cite{Hou}.
\end{rmk}

\bigskip
\noindent
{\bf Acknowledgements. } This work was financially supported by National
Natural Science Foundation of China (No.11301144, 11771122, 11801141).

 \end{document}